\documentclass[12pt]{article}

\usepackage[leqno]{amsmath}
\usepackage{amsfonts,amsthm,amssymb}
\usepackage[all]{xy}

\theoremstyle{plain}    
\newtheorem{theorem}{Theorem}[section]
\newtheorem{lemma}[theorem]{Lemma}
\newtheorem{proposition}[theorem]{Proposition}
\newtheorem{corollary}[theorem]{Corollary}
\newtheorem*{thm1}{Theorem~\ref{the: Convergence of EMSS over K(P,1)}}
\newtheorem*{thm2}{Theorem~\ref{the: Convergence of EMSS to fixed points over finite nilpotent space}}
\newtheorem*{IndependenceLem}{Lemma~\ref{lem: Independence of Base}}

\theoremstyle{definition}   
\newtheorem{definition}[theorem]{Definition}
\newtheorem{example}[theorem]{Example}
\newtheorem{remark}[theorem]{Remark}

\theoremstyle{remark}   

\newcommand{\Int}{\mathbb{Z}}   
\newcommand{\sphere}{\mathbb{S}}
\newcommand{\field}{\mathbb{F}} 

\newcommand{\cc}{\mathcal{C}}   
\newcommand{\Derived}{\mathbf{D}}   

\newcommand{\Hom}{\text{Hom}}   
\newcommand{\ext}{\text{Ext}}   
\newcommand{\tor}{\text{Tor}}   
\newcommand{\End}{\text{End}}   
\newcommand{\cone}{\mathbf{Cone}}  
\newcommand{\chains}{\mathrm{C}}    
\newcommand{\hocolim}{\mathrm{hocolim}} 
\newcommand{\colim}{\mathrm{colim}} 
\newcommand{\cell}{\mathrm{Cell}}   
\newcommand{\kangr}{\Omega} 
\newcommand{\B}{\mathrm{B}} 
\newcommand{\W}{\mathrm{W}} 

\newcommand{\stam}[1]{}         

\newcommand{\problem}[1]{ }
\newcommand{\elaboration}[1]{}

\begin{document}    

\title{Cellular approximations and the Eilenberg-Moore spectral-sequence}   
\author{Shoham Shamir}  
\date{\today}                   
\maketitle                      

\begin{abstract}
We set up machinery for recognizing $k$-cellular modules and $k$-cellular approximations, where $k$ is an
$R$-module and $R$ is either a ring or a ring-spectrum. Using this machinery we can
identify the target of the Eilenberg-Moore cohomology spectral sequence for a fibration in various cases.
In this manner we get new proofs for known results concerning the Eilenberg-Moore spectral sequence and generalize 
another result.

\end{abstract}

\section{Introduction}
\label{sec: Introduction}

Given a pointed topological space $A$, the $A$-cellularization of another pointed space $X$ is the closest
approximation of $X$ built from $A$ using pointed homotopy colimits.
This notion of $A$-cellularization was first treated systematically by Dror-Farjoun \cite{Farjoun} and
Chach\'{o}lski \cite{Chacholski}. Also fundamental to this subject is the work of Hirschhorn \cite{Hirschhorn}.
Dwyer and Greenlees translated this concept to an algebraic setting
in~\cite{Dwyer Greenlees}. Further work on cellular approximation was done by Dwyer,Greenlees and Iyengar in 
\cite{Dwyer Greenlees Iyengar}, here in the wider setting of $\sphere$-algebras, where $\sphere$ is the sphere
spectrum.

The algebraic setting of Dwyer and Greenlees from~\cite{Dwyer Greenlees} is as follows.
Let $R$ be a ring (with a unit, not necessarily commutative) and let $k$ and $C$ be $\Int$-graded chain-complexes
of left $R$-modules. We say $C$ is \emph{$k$-cellular} over $R$, if $\ext_R^*(k,N)=0$ implies $\ext_R^*(C,N)=0$ for 
every $R$-chain complex $N$. A $k$-cellular approximation of a complex $X$ is a map of $R$-complexes: 
$\cell_k^R X \to X$ such that $\cell^R_k X$ is the closest approximation of $X$ by a $k$-cellular complex
in the sense that the induced map: $\ext_R^*(k,\cell_k^R X) \to \ext_R^*(k,X)$ is an isomorphism.

The concept of $k$-cellular approximation extends effortlessly to the stable homotopy category, which is the
setting used here. Specifically we let $R$ be an $\sphere$-algebra in the sense of \cite{EKMM}, where $\sphere$ denotes the sphere-spectrum. Let $k$ and $C$ be left $R$-modules. Then the definitions above immediately apply also 
to the category of left $R$-modules. 


Dwyer, Greenlees and Iyengar gave in~\cite{Dwyer Greenlees Iyengar}, a formula for $k$-cellular
approximation. This formula holds when $k$ is a \emph{proxy-small} $R$-module 
(see definition~\ref{def: Proxy Small}), in
such a case the $k$-cellular approximation of an $R$-module $X$ is given by the following natural map 
\cite[theorem 4.10]{Dwyer Greenlees Iyengar}:
\[ \Hom_R(k,X) \otimes_{\End_R(k)} k \to X\]
The symbol $\otimes_R$ denotes the smash product of $R$-modules and $\Hom_R(-,-)$ stands 
for the function spectrum of $R$-modules. These functors are taken in the derived sense, i.e. we always 
take appropriate resolutions before applying them.

Dwyer, Greenlees and Iyengar's formula for $k$-cellular approximations can be used to understand the target
of the Eilenberg-Moore cohomology spectral sequence for a fibration, as we will now demonstrate.
Suppose $k$ is a commutative $\sphere$-algebra and let $F \to E \to B$ be fibration of spaces where $E$ and $B$ 
are connected. We make the following standard replacements. First, one can take replace $\kangr B$ by an equivalent
topological group, so we will assume $\kangr B$ is a topological group. Second, standard constructions show we can 
replace $F$ by an equivalent $\kangr B$-space (by passing, perhaps, to an equivalent fibration) so we will assume 
that $F$ itself is a right $\kangr B$-space. Let $\chains^*(F;k)$ denote the function spectrum 
$F_\sphere(\Sigma^\infty F_+,k)$ and let $k[\kangr B]$ denote the $\sphere$-algebra 
$k \wedge_\sphere \Sigma^\infty \kangr B_+$. Thus $\chains^*(F;k)$ is a left $k[\kangr B]$-module.
The same arguments of Dwyer and Wilkerson from \cite[2.10]{Dwyer Wilkerson} can be used to show there is a 
weak equivalence of $k$-modules:
\[ \Hom_{k[\kangr B]} (k,\chains^*(F;k)) \otimes_{\End_{k[\kangr B]}(k)} k \simeq
\chains^*(E;k) \otimes_{\chains^*(B;k)} k\]

On the other hand, recall the Eilenberg-Moore cohomology spectral sequence for this fibration (with coefficients in 
$k$), has the form:
\[ E^2_{p,q}= \tor^{H^*(B;k)}_{p,q}(H^*(E;k),k) \ \Rightarrow \ 
\pi_{p+q}(\chains^*(E;k) \otimes_{\chains^*(B;k)} k)\]
Thus, whenever $k$ is a proxy-small $k[\kangr B]$-module, the Eilenberg-Moore cohomology spectral sequence
computes the homotopy groups of the $k$-cellular approximation of $\chains^*(F;k)$ over $k[\kangr B]$.
In particular, if $\chains^*(F;k)$ is a $k$-cellular $k[\kangr B]$-module, then the target of the Eilenberg-Moore
cohomology spectral sequence is the cohomology of the $F$ with coefficients in $k$.

Using machinery for recognizing $k$-cellular modules and $k$-cellular approximations we get the following two
theorems. 

\begin{thm1} 
Fix a prime $p$. Let $N$ be a finite nilpotent group and let $P\subseteq N$ be the $p$-Sylow subgroup of $N$,
so $N \cong P \times H$ with the order of $H$ being prime to $p$. 
Let $F \to E \to \B N$ be a homotopy fibration sequence over the classifying space of $N$,
with $E$ being a connected space.
If $k$ is any commutative $\sphere$-algebra such that $\pi_0(k)$ is an $\field_p$-algebra, then:
\[ \chains^*(E;k) \otimes_{\chains^*(B;k)} k \simeq \chains^*(F_{h(H)})\]
Where $F_{h(H)}$ is the homotopy orbit space of $F$ with respect to the $H$-action on $F$.
In particular, if $N=P$ then \[ \chains^*(E;k) \otimes_{\chains^*(B;k)} k \simeq \chains^*(F;k) \]
\end{thm1} 

\begin{thm2} 
Fix a prime $p$. Let $B$ be a finite connected nilpotent space with a finite fundamental group $N=\pi_1(B)$.
Let $P\subseteq N$ be the $p$-Sylow subgroup of $N$, so $N \cong P \times H$.
Let $F \to E \to B$ be a homotopy fibration sequence over $B$, where $E$ is a connected space. Then:
\[ \tor_{-n}^{\chains^*(B;\field_p)}(\chains^*(E;\field_p), \field_p) = H^n(F;\field_p)^{H}\]
Where $H^n(F;\field_p)^{H}$ are the fixed points of the $H$-action on $H^n(F;\field_p)$.
\end{thm2} 
The first theorem deals with general multiplicative cohomology theories in characteristic $p$ and thus 
generalizes a result of Kriz from~\cite{Kriz}, showing convergence of the Eilenberg-Moore 
cohomology spectral sequence for a fibration over $\B (\Int/p)$, with coefficients in Morava $K$-theories.
The second result concerns only mod-$p$ cohomology. 
This result is hardly surprising as it is dual to Dywer's exotic convergence \cite{Dwyer},
which concerns convergence of the Eilenberg-Moore homology spectral sequence.
However, the proof we give here uses only cellularity arguments. We also give a different proof to a weaker
version of Dwyer's strong convergence result from~\cite{DwyerStrong}, again using only cellularity
arguments. This is done in proposition~\ref{pro: Different proof of Dwyer's strong convergence}.

Another result is a spectral sequence.
Assuming that $k$ is a commutative ring, we describe, in section~\ref{sec: Applications}, a spectral sequence with 
\[E^1_{p,q} = \pi_{2p+q}(\cell^{k[\pi_1 B]}_k H^p(F;k))\]
that converges to $\pi_{p+q}( \cell_k^{k[\kangr B]} \chains^*(F;k))$ when $B$ is a finite nilpotent space.
We demonstrate the use of this spectral sequence in lemma~\ref{lem: exact sequence for fibration with cyclic pi1}
and corollary~\ref{cor: when pi1B is cyclic of order 2}.

The bulk of this paper is concerned with setting up the necessary machinery for recognizing $k$-cellular 
modules and $k$-cellular approximations. This machinery is interesting in itself and we mention one basic
result. Let $R \to S$ be a map of $\sphere$-algebras and let $k$ be an $S$-module.
The independence of base lemma (lemma~\ref{lem: Independence of Base}) relates $k$-cellular approximation of 
$S$-modules with their $k$-cellular approximation as $R$-modules.
\begin{IndependenceLem} 
Let $R \to S$ be a map of $\sphere$-algebras.
\begin{enumerate} 
\item (\emph{Strong independence of base})
Let $k$ be an $S$-module. If $S \otimes_R k$ is $k$-cellular as an $S$-module, then for any $S$-module $X$
the map \mbox{$\cell_k^S X \to X$} is a $k$-cellular approximation of $X$ as an $R$-module.
In particular $\cell_k^R X \simeq \cell_k^S X$ as $R$-modules.
\item (\emph{Weak independence of base})
Let $k$ be an $R$-module. If $S \otimes_R k$ is $k$-cellular as an $R$-module, then for any $S$-module $X$ 
the map $\cell_{S \otimes_R k}^S X \to X$ is a $k$-cellular approximation of $X$ as an $R$-module.
In particular $\cell_k^R X \simeq \cell_{S \otimes_R k}^S X$ as $R$-modules.
\end{enumerate}
\end{IndependenceLem} 
This lemma is the basis for results on $k$-cellular approximation over a group-ring $kG$, where $G$
is a simplicial, or discrete, group and $kG$ denotes the $\sphere$-algebra $\chains_*(G;k)=k \wedge_\sphere 
\Sigma^\infty G_+$. See for example proposition~\ref{pro: cellularization over a group-ring}.

\subsection{Layout of this Work}
We begin by introducing the notions of cellularity and cellular approximation in 
section~\ref{sec: Cellular Approximations}. There we also recall some basic properties of cellular approximations 
and give elementary examples.

The necessary machinery for recognizing $k$-cellular modules and building $k$-cellular approximations
is set up is sections~\ref{sec:Change of Rings} and \ref{sec:Cellular Approximation in Nilpotent Groups}.
In section~\ref{sec:Change of Rings} we prove the independence of base lemma and
apply it to get results concerning $k$-cellular approximation over a group ring $kG$.
Section~\ref{sec:Cellular Approximation in Nilpotent Groups} deals with $k$-cellular approximation
over finite $p$-groups and finite nilpotent groups.

Finally, in section~\ref{sec:Fibrations}, we come to the Eilenberg-Moore spectral sequence, applying the machinery
we constructed earlier to identify it's target. The last section (\ref{sec: Applications}) is devoted to constructing
the spectral sequence and demonstrating its use.

\subsection{Setting and Conventions}
\label{sub:Setting and Conventions}
We work in the category of left $R$-modules where $R$ is an $\sphere$-algebra in the sense of~\cite{EKMM} and
$\sphere$ stands for the sphere spectrum. However, in some instances we will prefer to use an ordinary ring $R$
and work in the category of differential graded left $R$-modules (i.e. $\Int$-graded chain complexes of left 
$R$-modules). We justify this in the following way.
Recall that for a ring $R$, the $\sphere$-algebra $HR$ is the appropriate Eilenberg-Mac Lane $\sphere$-algebra.
A result of Schwede and Shipley \cite[5.1.6]{SchwedeShipley}, shows the
model category of $HR$-modules is Quillen equivalent to the projective model category of unbounded chain
complexes over $R$. In this way one can pass between $HR$-modules and differential graded $R$-modules
without fear and we will usually not distinguish between the two.
To keep the terminology consistent, when $R$ is a ring the term $R$-module will mean a differential
graded $R$-module and we will refer to the $R$-modules in the classical sense as \emph{discrete $R$-modules}
(this is the same terminology as in \cite{Dwyer Greenlees Iyengar}).

We follow~\cite{Dwyer Greenlees Iyengar} in both terminology and notation.
Thus, as in~\cite{Dwyer Greenlees Iyengar}, we use the notation $-\otimes_R-$ for the smash product of $R$-modules.
Similarly we use $\Hom_R(-,-)$ to denote the function spectrum of $R$-modules.
All functors in this paper are taken in the derived sense, in particular $-\otimes_R-$ and $\Hom_R(-,-)$.
This means that we always assume to have replaced our modules by appropriate resolutions before applying
the functor in question.

Given $R$-modules $k$ and $M$, we denote the $k$-cellular approximation of an $R$-module $M$ by $\cell_k^R M$,
or by $\cell_k M$ whenever $R$ is clear from the context.

The homotopy groups of an $R$-module $M$ are the stable homotopy groups of its underlying spectrum, and are
denoted as usual by $\pi_i(M)$. 
The derived (or homotopy) category of $R$-modules is denoted by $\Derived_R$.
Since this is a triangulated category, we sometimes use the term \emph{triangle} to indicate a homotopy cofibration
sequence. If $f:M \to N$ is a map of $R$-modules, then $\cone (f)$ denotes the mapping cone of $f$.
We say two $R$-modules are \emph{equivalent} if they are isomorphic in $\Derived_R$, and use the symbol
$\simeq$ to denote this. Similarly two maps between $R$-modules are said to be equivalent if they become
equal when passing to the derived category $\Derived_R$.
An $R$-module $M$ is called \emph{connective} if $\pi_i(M)=0$ for all $i<0$. We say that $M$ is
\emph{bounded above} (respectively \emph{bounded below}) if there exists some index $n$ such that $\pi_i(M)=0$
for all $i>n$ (respectively $i<n$).

If $k$ is a commutative $\sphere$-algebra and $X$ is a based space, then
the \emph{chains of $X$} with coefficients in $k$ is the spectrum:
\[ \chains_*(X;k) = k \otimes_\sphere \Sigma^\infty X\]
When $k$ is understood from the context we will simply use the notation $\chains_*(X)$.
If $X$ has no base point we add a disjoint base point to $X$ before taking the chains, i.e.
$\chains_*(X;k) = k \otimes_\sphere \Sigma^\infty X_+$.
Similarly, the \emph{cochains} of based space $X$ with coefficients in $k$ is the function spectrum:
\[ \chains^*(X;k) = \Hom_\sphere(\Sigma^\infty X,k)\]
All our spaces are assumed to have the homotopy type of CW-spaces. We use $pt$ to denote the space with a single 
point. For a based space $X$, the space $\kangr X$ is a topological group weakly equivalent to the loop-space
of $X$. 

In many cases the ring $R$ we work over will be a group-ring. If $G$ is a topological group 
and $k$ is a commutative $\sphere$-algebra, the group-ring $kG$ is $\chains_*(G;k)$,
which is an $\sphere$-algebra.
If $k$ is a commutative ring and $G$ a discrete group, then $Hk[G]$ is equivalent to the Eilenberg-Mac Lane spectrum
of the usual group-ring $kG$. When $R$ is a group-ring $kG$, the notation $\cell^G(-)$ will stand for 
$\cell^{kG}_k(-)$, so long as the ground $S$-algebra $k$ is understood. 
Note that a map $f:X \to Y$ of $G$-spaces which yields an isomorphism on all homotopy groups will induce a weak 
equivalence $\chains_*(f;k): \chains_*(X;k) \to \chains_*(Y;k)$ of $kG$-modules.

The classifying space of a topological group $G$ is denoted in the usual manner by $\B G$ and
$\W G$ is a contractible free $G$-space, such as the space described by the Rothenberg and Steenrod 
\cite{RothenbergSteenrod}.

\subsection{Acknowledgements}
I am grateful to Prof. Emmanuel Dror-Farjoun for his guidance, support and for many useful suggestions. 
Most notably, it was he
who led me to consider the Eilenberg-Moore spectral sequence in the context of cellular approximations.
Part of this work was done while visiting the Centre De Recerca M\'{a}tematica in Barcelona under a Marie Curie
scholarship. I would like to thank Prof. Carles Casacuberta who enabled this visit. 

\section{Cellular Approximations}
\label{sec: Cellular Approximations}
Throughout this section we fix an $\sphere$-algebra $R$ and an $R$-module $k$.
In this section we give the necessary background on cellular approximations and some of their basic properties.
After defining $k$-cellulatrity and $k$ cellular-approximations, we recall from \cite[4.6]{Dwyer Greenlees Iyengar} 
the definition of proxy-smallness. By imposing the condition of proxy-smallness on $k$ one gets better handle on 
$k$-cellular approximations, as shown by Dwyer, Greenlees and Iyengar in that paper
\cite[4.10]{Dwyer Greenlees Iyengar}.
We end this section with several examples in which $k$-cellularity has a simple description.

We recall the following definitions from~\cite{Dwyer Greenlees Iyengar}.
\begin{definition}
\label{def: cellular approximations}
A map $U \to V$ of $R$-modules is called a \emph{$k$-equivalence} if the induced map of $\sphere$-modules:
$\Hom_R(k,U) \to \Hom_R(k,V)$ is an equivalence.
An $R$-module $N$ is called \emph{$k$-null} if $\Hom_R(k,N)\simeq 0$.
An $R$-module $M$ is called \emph{$k$-cellular} if, for every $k$-null module $N$, $\Hom_R(M,N) \simeq 0$.
Equivalently, $M$ is \emph{$k$-cellular} if for every $k$-equivalence $f:U \to V$ the induced map
$\Hom_R(M,U) \to \Hom_R(M,V)$ is an equivalence.

If $C \to X$ is a $k$-equivalence and $C$ is $k$-cellular then $C$ is called a \emph{$k$-cellular approximation
of $X$} (or a \emph{$k$-cellularization of $X$}). The module $C$ is denoted by $\cell^R_k (X)$ or simply
$\cell_k (X)$ if $R$ is clear from the context.
\end{definition}

\begin{remark}
Note that the property of being a $k$-cellular approximation can be defined solely in the derived category 
$\Derived_R$. Also, from definition~\ref{def: cellular approximations} it is easy to see that a $k$-cellular 
approximation of any given module is unique up to a unique isomorphism in the derived category $\Derived_R$. 
\end{remark}

An important example of cellular approximation was given by Dwyer and Greenlees in~\cite{Dwyer Greenlees}.
Suppose $R$ is a Noetherian commutative ring, let $I$ be an ideal of $R$ and take $k$ to be the
ring $R/I$. Dwyer and Greenlees showed that for any $R$-module $X$ the homotolgy groups of $\cell_k X$ 
are isomorphic to the local cohomology of $X$ at $I$, i.e.:
\[ \pi_{-n}(\cell_k X) \cong H_I^n (X) = \pi_{-n}(\colim_m \Hom(R/I^m,X)) \]

Here is a simpler example. Suppose $R$ is a ring (not necessarily commutative) and $I$ is an ideal of $R$.
Recall that a discrete $R$-module $M$ is said to be $I$-nilpotent if there is a filtration
$0=M_0 \subset M_1 \subset \cdots \subset M_n=M$ such that $M_i/M_{i+1}$ is an $R/I$-module. It is a simple
matter to show that every $I$-nilpotent module is $R/I$-cellular.

There is a well known, equivalent definition of $k$-cellularity, which we will use extensively.
To state this definition we must first recall what are thick and localizing subcategories.
\begin{definition}
Let $\cc$ be a full subcategory of the category of $R$-modules.
The subcategory $\cc$ is called \emph{thick} if it is closed under equivalences, 
triangles and retracts (direct summands).
Closure under triangles means that given a triangle $X \to Y \to Z$ where two of the modules belong to $\cc$, 
then so does the third. In particular a thick subcategory is closed under suspensions and finite
coproducts. A thick subcategory is called \emph{localizing} if in addition it is closed under arbitrary coproducts.

Let $A$ be an $R$-module. The smallest thick subcategory containing $A$ is called the \emph{thick subcategory
generated by $A$}. The \emph{localizing subcategory generated by $A$} is similarly defined. 
Following \cite{Dwyer Greenlees Iyengar}, we say an $R$-module $B$ is \emph{finitely built} by $A$ if
$B$ belongs to the thick subcategory generated by $A$. An $R$-module which is finitely built by $R$
 is also called a \emph{small} $R$-module. Finally, as in \cite{Dwyer Greenlees Iyengar},
we say $B$ is \emph{built} by $A$ if $B$ belongs to the localizing subcategory generated by $A$.
Clearly, if $B$ is finitely built from $A$ then $B$ is also built from $A$.
\end{definition}

We can now state the equivalent definition for $k$-cellular modules:
an $R$-module is $k$-cellular if and only if it is built from $k$. One part of this equivalence is easy to see,
namely that any  $R$-module built from $k$ is $k$-cellular. Using the techniques of \cite{Farjoun} or \cite[5.1.5]{Hirschhorn}, one can show that for any $R$-module $X$ there exists a $k$-cellular approximation. This 
fact proves the other direction. 

\begin{remark}
\label{rem: cellularization of a set of objects}
Let $S=\{A_i\}_{i\in I}$ be a set of $R$-modules. Say an $R$-module $M$ is \emph{built by
$S$} (equiv. \emph{$S$-cellular}) if $M$ belongs to the localizing subcategory generated by $S$.
Clearly, this subcategory is equivalent to the localizing subcategory generated by the
$R$-module $A=\oplus_{i\in I}A_i$. Therefore being $A$-cellular is the same as being $S$-cellular.
\end{remark}

As noted above, given $R$-modules $k$ and $X$ one can always construct a $k$-cellular approximation for $X$.
In general these constructions involve transfinite induction and a small object argument.
To get a more manageable construction of $k$-cellular approximation, it is necessary to impose conditions on $k$.
One such condition, given by Dwyer, Greenlees and Iyengar in~\cite[4.6]{Dwyer Greenlees Iyengar} is
that $k$ be proxy-small. We recall the definition of proxy-smallness in~\ref{def: Proxy Small} below.
Dwyer, Greenlees and Iyengar proved in~\cite[4.10]{Dwyer Greenlees Iyengar} that when $k$ is a proxy-small $R$-module the $k$-cellular approximation of any $R$-module $X$ is given by the natural map: 
\[ \Hom_R(k,X) \otimes_{\End_R(k)} k \to X \]

\begin{definition}
\label{def: Proxy Small}
Call an $R$-module $k$ \emph{proxy-small}
if there exists an $R$-module $K$ such that $K$ is finitely built
from $R$, $K$ is also finitely built from $k$ and $K$ builds $k$. The module $K$ is
called a \emph{Koszul complex associated to $k$} 
(this is the term used in~\cite{Dwyer Greenlees Iyengar}, note that in
\cite{Dwyer Greenlees Iyengar Finiteness} the module $K$ is called \emph{a witness that $k$ is proxy-small}).
\end{definition}
\begin{remark}
Suppose $k$ is proxy-small with an associated Koszul complex $K$.
Because $K$ and $k$ build each other, an $R$-module $M$ is built from $k$ if and only if $M$ is built from $K$.
Moreover, a map of $R$-modules is a $K$-equivalence if and only if it is a $k$-equivalence.
Therefore a $K$-cellular approximation is the same as a $k$-cellular approximation.
\end{remark}

In \cite[example 5.9]{Dwyer Greenlees Iyengar}, Dwyer, Greenlees and Iyengar show that for any finite group
$G$ and any commutative ring $k$, the trivial $kG$ module $k$ is proxy-small. 
In proposition~\ref{pro: k proxy-small over kN} below we will see that if $G$ is a finite nilpotent group and $k$ is
a commutative $\sphere$-algebra such that $\pi_0(k)$ is an $\field_p$-algebra, then $k$ is a proxy-small $kG$-module.

We end this section by showing two cases where $k$-cellularity is relatively simple.
In both cases $k$ is an $\sphere$-algebra which is given an $R$-module structure via a map of 
$\sphere$-algebras $R \to k$.

In the first case $R$ is a connective $\sphere$-algebra, i.e. $\pi_n(R)=0$ for all $n<0$. 
Recall that for an abelian group $A$, the appropriate Eilenberg-Mac Lane $\sphere$-module is denoted $HA$.
By~\cite[IV.3.1]{EKMM}, there exists map of $\sphere$-algebras $R \to H\pi_0(R)$
which yields the obvious isomorphism on the $\pi_0$-level. In particular, $H\pi_0(R)$ is an $R$-module. We show
that when $R$ is connective, every bounded above $R$-module is $H\pi_0(R)$-cellular.
This result is dual to a result of A.K. Bousfield from~\cite[lemma 3.3]{Bousfield spectra} over $\sphere$.
\begin{proposition}
\label{pro: coconnective module over a connective algebra A is A_0-cellular}
Let $R$ be a connective $\sphere$-algebra, then every bounded above $R$-module is $H\pi_0(R)$-cellular.
\end{proposition}
\begin{proof}
Let us denote the $\sphere$-algebra $H\pi_0(R)$ by $k$. Note that for any $R$-module $M$, the modules
$H \pi_i(M)$ are $k$-modules and hence also $R$-modules. Let $X$ be a bounded above $R$-module.
Connectivity of $R$ implies the existence of Postnikov sections in the category of $R$-modules
(see~\cite[lemma 3.2]{Dwyer Greenlees Iyengar}).
This allows us to use an argument of Dwyer and Greenlees 
from~\cite[proof of proposition 5.2]{Dwyer Greenlees}, showing that $X$ is built from the $R$-modules 
$\{H \pi_n(X)\}_{n \in \Int}$. We detail this argument below.

Following Dwyer and Greenlees \cite{Dwyer Greenlees}, we 
denote by $M\langle -\infty,j\rangle $ the $j$'th Postnikov section of an $R$-module $M$.
Thus, there is a natural map of $R$-modules $p_j:M \to M\langle -\infty,j\rangle$ such that
$\pi_k(p_j)$ is an isomorphism for $k\leq j$ and $\pi_m(M\langle -\infty,j\rangle)$ is zero for $m>j$.
Define $M\langle i,\infty\rangle $ to be the homotopy fiber of the map 
$M \to M\langle -\infty,i-1\rangle $, yielding a triangle:
\[ M\langle i,\infty\rangle  \xrightarrow{c_i} M \xrightarrow{p_{i-1}} M\langle -\infty,i-1\rangle \]
Define $M\langle i,j\rangle $ (where $i\leq j$) to be \[(M\langle i,\infty\rangle )\langle -\infty,j\rangle\] 
(this notation agrees with the notation of \cite[proof of proposition 5.2]{Dwyer Greenlees} for $R$-chain complexes,
up to equivalence). Particularly useful is the triangle:
\[ M\langle i,j\rangle  \to M\langle i-1,j\rangle  \to M\langle i-1,i-1\rangle  \]
Also note that $M\langle i,i\rangle  \simeq H\pi_i(M)$.

Here is the argument from~\cite[proof of proposition 5.2]{Dwyer Greenlees}.
Since $X$ is bounded above, there is an index $j$ such that $\pi_n(X)=0$ for all $n > j$. 
First, it will be proved that for any $i \leq j$, the $R$-module $X \langle i, j \rangle$ is built from 
$\{H \pi_n(X)\}_{i\leq n\leq j}$. This statement is clearly true for $i=j$. 
It is true for $i<j$ by induction on $i$, using the triangle:
\[ X\langle i,j\rangle  \to X\langle i-1,j\rangle  \to X\langle i-1,i-1\rangle \] 
Second, there is an obvious equivalence:
$X\langle i,j\rangle  \xrightarrow{\simeq} X\langle i,\infty\rangle $ and one can form a telescope:
\[ X\langle j,j\rangle  \to X\langle j-1,j\rangle   \to \cdots \to X\langle i,j\rangle  \to \cdots \]
with obvious maps from this telescope to $X$. It is easy to see that the map
\[\hocolim_{(i\to -\infty)} X\langle i,j\rangle  \to X\]
is an equivalence. Since $X\langle i,j\rangle $ is built by $\{H \pi_n(X)\}_{i\leq n\leq j}$, we conclude
$X$ is built by $\{H \pi_n(X)\}_{ n\leq j}$.

To finish the proof, note that for any $n$, $H \pi_n(X)$ is a $k$-module and therefore it is clearly built by
$k$ in the category of $k$-modules. Hence $H \pi_n(X)$ is also built by $k$ in the category of $R$-modules.
\end{proof}

\begin{example}
Let $k$ be a commutative ring (in the classical sense), let $G$ be a topological group and let $X$ be a $G$-space.
Denote by $G_0$ the discrete group $\pi_0(G)$. It is easy to see that $\pi_0(kG) \cong k G_0$ and that
the map of $\sphere$-algebras $kG \to kG_0$ is induced by the map of topological groups $G \to \pi_0(G)$.
The group-ring $kG$ is connective, the $kG$-module $\chains^*(X;k)$ is coconnective, therefore,
by proposition~\ref{pro: coconnective module over a connective algebra A is A_0-cellular} above,
$\chains^*(X;k)$ is $kG_0$-cellular. In fact, $\chains^*(X;k)$ is built
from the cohomology groups $H^n(X;k)$ which are $kG_0$-cellular as $kG$-modules 
(because $H^n(X;k)$ is a $kG_0$-module).
\end{example}

\begin{example}
\label{exa:For a finite nilpotent space tildeB is a k-cellular pi_1(B)-space}
Let $B$ be a finite connected nilpotent space with a fundamental group $G=\pi_1(B)$.
A result of Hilton, Mislin and Roitberg \cite[II.2.18]{Hilton Mislin Roitberg} shows $G$ operates
nilpotently on the modules $H_*(\tilde{B})$, where $\tilde{B}$ is the universal cover of $B$.
Since $B$ is a finite space we see $\chains_*(\tilde{B};\Int)$ is a bounded above $\Int G$-module 
(i.e. a bounded above $\Int G$-chain complex).
From the proof of proposition~\ref{pro: coconnective module over a connective algebra A is A_0-cellular} we see that
$\chains_*(\tilde{B};\Int)$
is built by it's homology groups. Every $H_n(\tilde{B})$ is a nilpotent $\Int G$-module and hence $\Int$-cellular.
Therefore $\chains_*(\tilde{B};\Int)$ is a $\Int$-cellular $\Int G$-chain complex.

This generalizes easily for any commutative ring $k$, by the following observation.
It is a simple matter to show that applying
the functor $k \otimes_\Int -$ turns $\Int$-cellular $\Int G$-modules into $k$-cellular $kG$-modules.
Since $\chains_*(\tilde{B};k)$ is equivalent to $k \otimes_\Int \chains_*(\tilde{B};\Int)$,
it is a $k$-cellular $kG$-chain complex.

Moreover, this result can be generalized for any $\sphere$-algebra $k$ that is a commutative, 
connective and bounded above, by the following argument.
Using the Atiyah-Hirzebruch spectral sequence (see e.g.~\cite[IV.3.7]{EKMM}) we see that 
$\pi_n(\chains_*(\tilde{B};k))$
are nilpotent $\pi_0(kG)$-modules. Since $k$ is bounded above, $\chains_*(\tilde{B};k)$ is also bounded above.
Therefore $\chains_*(\tilde{B};k)$ is built from the nilpotent $\pi_0(kG)$-modules
$H \pi_n(\chains_*(\tilde{B};k))$. Since the modules $\pi_n(\chains_*(\tilde{B};k))$ 
are nilpotent, they are built by $\pi_0(k)$ in the category 
of $\pi_0(kG)$-modules. This implies $H \pi_n(\chains_*(\tilde{B};k))$ is built by $H\pi_0(k)$ in the category of 
$kG$-modules. Since $H\pi_0(k)$ is a $k$-module, it is built by $k$ (also over $kG$).
Therefore $\chains_*(\tilde{B};k)$ is a $k$-cellular $kG$-module.
\end{example}

\begin{example}
Here is a simple example of a connective $\sphere$-algebra $R$ and an unbounded $R$-module $M$ which is not 
$H\pi_0(R)$-cellular. Let $k=H\pi_0(R)$ and suppose $R$ is such that
$\pi_*(R)$ is isomorphic to the graded ring $k_*[x]$, with $x$ in dimension greater than 0.
Let $X$ be the homotopy colimit of the telescope:
\[R \xrightarrow{x} \Sigma^{-|x|}R \xrightarrow{x} \Sigma^{-2|x|}R\xrightarrow{x} \cdots\]
where $R \xrightarrow{x} \Sigma^{-|x|}R$ is a map representing multiplication by $x$. 
Clearly, $X$ is not bounded above.
It is easy to show that $X$ is $k$-null and in particular $X$ is not $k$-cellular.
\elaboration{From \cite[proposition 3.9]{Dwyer Greenlees Iyengar} the mapping cone of the map 
$\Sigma^{|x|}R \xrightarrow{x} R$ is equivalent to $k$.}
\end{example}

We now turn to the second case where $k$-cellularity is particularly simple.
Here we suppose $k$ is a retract of $R$ in the derived category $\Derived_R$ and we obtain the following
simple result, which will be put to good use in section~\ref{sec:Cellular Approximation in Nilpotent Groups}.
\begin{lemma}
\label{lem: k is a retract of R}
Suppose $k$ is an $R$-module via an $\sphere$-algebra map $a:R \to k$ and suppose there exists in $\Derived_R$ a map 
$b:k \to R$ such that $a b = id_{k}$ (in $\Derived_R$). Then, for any $R$-module $X$, the $k$-cellular 
approximation of $X$ is given by the natural map:
\[Hom_R(k,X) \xrightarrow{a^*} X\]
\end{lemma}
\begin{remark}
\label{rem: cofibrant replacement of a bimodule}
Recall that $\Hom_R(k,X)$ is taken in the derived sense, i.e. we take a cofibrant replacement of $k$ before applying
the functor $\Hom_R(-,X)$. But now we must explain why this construction results in an $R$-module.
The point is that we can replace $k$ by an $R$-bimodule which is cofibrant as a left $R$-module.
So the remaining right action on this bimodule gives the desired left $R$-module structure on the function
complex $\Hom_R(k,X)$.
\end{remark}
\begin{proof}
We start by showing the map $a \otimes 1:k \to k \otimes_R k$ is an equivalence.
Clearly $(a\otimes 1) (b\otimes 1)= id_k$ in $\Derived_R$ i.e. $b \otimes 1$ is
a right inverse (in $\Derived_R$) for $a \otimes 1$. We need only show that $a \otimes 1$ has a left inverse
in $\Derived_R$, since this would imply $b \otimes 1$ is both a left and right inverse for $a \otimes 1$.

The map $a$ induces a functor $a^*:\Derived_k \to \Derived_R$ which is right adjoint to the functor 
$a_*=R/I \otimes_R -$. If $A$ is an $R$-module and $B$ is a $k$-module then the adjoint of a map
$f:k \otimes_R A \to B$ is the composition $(a^*f) \circ (a\otimes_R 1): A \to a^*B$.
Let $\mu:k \otimes_R k \to k$ be the map of $R/I$-modules that is adjoint to the identity
map $id_k$ via the adjunction above. Set $m=a^*\mu$, then the composition:
\[k \xrightarrow{a\otimes 1} k \otimes_R k \xrightarrow{m} k\]
is the map adjoint to $\mu$, hence it is the identity. This proves $m$ is a left inverse for $a \otimes 1$ 
in $\Derived_R$.

Given an $R$-module $X$ consider the map $a^*_X:\Hom_R(k,X) \to X$. Since $\Hom_R(k,X)$ is a
$k$-module, it is $k$-cellular. All that is left is to show $a^*_X$ is a $k$-equivalence.
Applying the functor $\Hom_R(k,-)$ to the map $a^*_X$ gives:
\[\Hom_R(k,a^*_X): \Hom_R(k,\Hom_R(k,X)) \to \Hom_R(k,X)\]
which is equivalent to the map:
\[ \Hom_R(k \otimes_R k ,X) \xrightarrow{(a \otimes 1)^*} \Hom_R(k,X)\]
Since $a \otimes 1$ is an equivalence, so is $(a \otimes 1)^*=\Hom_R(k,a \otimes 1)$.
\end{proof} 

\section{Change of Rings}
\label{sec:Change of Rings}
Suppose we are given a map of $\sphere$-algebras $R \to S$ and an $S$-module $k$. Given another $S$-module
$X$ there are two possible $k$-cellular approximations of $X$ we can consider:
the $k$-cellular approximation over $S$, which is $\cell_k^S(X)$ and the $k$-cellular
approximation over $R$, which is $\cell_k^R(X)$. We open this section by examining the relation between these
two cellular approximation in lemma~\ref{lem: Independence of Base} below.

\begin{lemma}
\label{lem: Independence of Base}
Let $R \to S$ be a map of $\sphere$-algebras.
\begin{enumerate}
\item (\emph{Strong independence of base})
Let $k$ be an $S$-module. If $S \otimes_R k$ is $k$-cellular as an $S$-module, then for any $S$-module $X$
the map \mbox{$\cell_k^S X \to X$} is a $k$-cellular approximation of $X$ as an $R$-module.
In particular $\cell_k^R X \simeq \cell_k^S X$ as $R$-modules.
\item (\emph{Weak independence of base})
Let $k$ be an $R$-module. If $S \otimes_R k$ is $k$-cellular as an $R$-module, then for any $S$-module $X$ 
the map $\cell_{S \otimes_R k}^S X \to X$ is a $k$-cellular approximation of $X$ as an $R$-module.
In particular $\cell_k^R X \simeq \cell_{S \otimes_R k}^S X$ as $R$-modules.
\end{enumerate}
\end{lemma}
\begin{remark}
We have mentioned before the relation between cellular approximation and local cohomology of commutative 
Noetherian rings, which was proved in~\cite{Dwyer Greenlees}.
The \emph{independence of base} property for local cohomology of commutative Noetherian rings
is the following:
Let $R \to S$ be a map of commutative Noetherian rings and let $I$ be an ideal of $R$. Then 
for every $S$-module $M$, the $I$-local cohomology groups of $M$, as an $R$-module, are isomorphic to the 
$IS$-local cohomology of groups $M$ as an $S$-module, where $IS$ is the ideal of $S$ generated by $I$. 

Using Dwyer and Greenlees' result \cite[proposition 6.10]{Dwyer Greenlees} one can
restate this property in cellular approximation terms. Indeed, it is easy to show that
this independence of base property is nothing other than an equivalence of the cellular approximations:
$\cell^R_{R/I} M \simeq \cell^S_{S/IS} M$ as $R$-modules.

We have dubbed the two parts of lemma~\ref{lem: Independence of Base} as strong and weak independence of base, 
because they describe properties are analogous to the independence of base property for local cohomology. In fact,
independence of base for local cohomology can be proved using the weak independence of base for cellular
approximation and Dwyer and Greenlees' result.
\end{remark}

\begin{proof}[Proof of lemma~\ref{lem: Independence of Base}]
We begin by making the following observation.
Suppose, as above, that $R \to S$ is a map of $\sphere$-algebras and $k$ is an $S$-module. Then every
$S$-module that is $k$-cellular over $S$ is also $k$-cellular over $R$. This follows easily from
the fact that the localizing subcategory of $S$-modules generated by $k$ is contained in the
localizing subcategory of $R$-modules generated by $k$.

We turn to prove the strong independence of base property. Suppose $S \otimes_R k$ is $k$-cellular as an 
$S$-module. Since $\cell_k^S X$ is $k$-cellular as an $S$-module, it is also $k$-cellular
as an $R$-module. We need only show that the map $\mu_X:\cell_k^S X \to X$ is a $k$-equivalence of $R$-modules.
Using adjunctions between the $\otimes$ and $\Hom$ functors (see~\cite[III.6.3]{EKMM}) we get:
\[ \Hom_R(k,\mu_X) \simeq \Hom_R(k,\Hom_S(S,\mu_X)) \simeq \Hom_S(S\otimes_R k, \mu_X)\]
The map $\mu_X$ is a $k$-equivalence of $S$-modules. The $S$-module $S \otimes_R k$ is $k$-cellular
and therefore $\mu_X$ is also an $S \otimes_R k$-equivalence. This proves $\Hom_R(k,\mu_X)$ is an equivalence.

We now prove the weak independence of base property. Let $k$ be an $R$-module and suppose
$S \otimes_R k$ is $k$-cellular as an $R$-module.
As above, $\cell_{S \otimes_R k}^S X$ is built by $S \otimes_R k$ also over $R$.
Since $S \otimes_R k$ is built by $k$ over $R$, we see $\cell_{S \otimes_R k}^S X$ is also built by $k$ over $R$.
We are left with showing the map $\nu_X:\cell_{S \otimes_R k}^S X \to X$ is 
a $k$-equivalence of $R$-modules. Using the same adjunctions as before, we get:
\[ \Hom_R(k,\nu_X) \simeq \Hom_S(S\otimes_R k, \nu_X)\]
Since $\Hom_S(S\otimes_R k, \nu_X)$ is an equivalence, we see that $\nu_X$ is a $k$-equivalence of $R$-modules.
\end{proof}

The main application of lemma~\ref{lem: Independence of Base} above is to group-rings.
We start with discrete groups.
Fix a commutative $\sphere$-algebra $k$ and a short exact sequence of discrete groups: $N \to G \to Q$.
The augmentation map of $\sphere$-algebras $kG \to k$ makes $k$ into a $kG$-module.
We denote $k$-cellular approximation over the group-ring $kG$ by $\cell^G$ whenever $k$ is clear from the
context. In similar fashion, we will refer to $kG$-modules as $G$-modules whenever $k$ is clear from the context.
We say $\cell^G$ is \emph{trivial} if for every $kG$-module $X$ the map $\cell^G X \to X$ is an equivalence.

The maps of groups induce maps $kN \to kG \to kQ \to k$ of $\sphere$-algebras. For any $G$-module
$X$ there are two possible $k$-cellular approximations: $\cell^G X$ and $\cell^N X$, each over a different ring.
Similarly, for a $Q$-module $Y$ we can consider two $k$-cellular approximations: $\cell^Q Y$ and $\cell^G Y$.
The relations between the various possible $k$-cellular approximations are given in 
proposition~\ref{pro: cellularization over a group-ring} below.

\begin{proposition}
\label{pro: cellularization over a group-ring}
Let $N \to G \to Q$ be a short exact sequence of discrete groups and let $k$ be a commutative $\sphere$-algebra.
\begin{enumerate}
\item 
For every $G$-module $X$ the map $\cell^G_{kQ}X \to X$ is a $k$-cellular approximation of 
$X$ as an $N$-module. In particular there is an weak equivalence of $N$-modules:
$\cell^N X \simeq \cell^G_{kQ}X$.
\item 
If $\cell^Q$ is trivial then for every $G$-module $X$ the map $\cell^G X \to X$ is a 
$k$-cellular approximation of $X$ as an $N$-module. In particular there is an equivalence of $N$-modules:
$\cell^G X \simeq \cell^N X$. 
\item
If both $\cell^Q$ and $\cell^N$ are trivial, then $\cell^G$ is also trivial.
\item
If $\chains_*(\B N;k)$ is a $k$-cellular $Q$-module, then for any $Q$-module $Y$ the map $\cell^Q Y \to Y$ is a 
$k$-cellular approximation of $Y$ as a $G$-module. In particular there is an equivalence
of $G$-modules: $\cell^Q Y \simeq \cell^G Y$.
\end{enumerate}
\end{proposition}

Before proving this proposition we need the following lemma.

\begin{lemma}
The module $kG \otimes_{kN} k$ is $k$-cellular as an $N$-module.
\end{lemma}
\begin{proof}
Clearly $G$ is isomorphic to $\coprod_{q\in Q} N$ as an $N$-space. Hence $KG \simeq \coprod_{q\in Q} kN$ as
a $kN$-module. The result follows.
\end{proof}

\begin{proof}[Proof of proposition~\ref{pro: cellularization over a group-ring}]
The previous lemma implies that for the map $kN \to kG$ and the $G$-module $k$ the
weak independence of base property holds (lemma~\ref{lem: Independence of Base}). This proves the first item
on our list.

To prove the second item we suppose $\cell^Q$ is trivial. This implies $kQ$ is $k$-cellular as a $Q$-module.
Since $k$ is clearly $Q$-cellular as a $Q$-module, we see $k$ and $kQ$ build each other over $kQ$.
But then $k$ and $kQ$ build each other also over $kG$.
We see that for every $G$-module $X$ the map $\cell^G_{kQ}X \to X$ is a $k$-cellular approximation of 
$X$ as a $G$-module.
To finish the proof of the second part of the proposition we apply the first part that, which was proved above.

Now suppose both $\cell^Q$ and $\cell^N$ are trivial and let $X$ be a $G$-module.
Because $\cell^Q$ is trivial, the natural map $\cell^G X \to X$ of $G$-modules is a
$k$-cellular approximation of $X$ as an $N$-module. Since $\cell^N$ is trivial, this map is an equivalence of
$N$-modules, hence it is also an equivalence of $G$-modules. This proves the third item on our list.

For the last item note that $\chains_*(\B N;k) \simeq kQ \otimes_G k$. Thus the result follows from 
strong independence of base (lemma~\ref{lem: Independence of Base})
\end{proof}

\begin{example}
Let $N \to G \to Q$ be a short exact sequence of finite groups with $N$ being a central subgroup of $G$
and $Q$ an abelian group. Let $k$ be a commutative ring. It is easy to show that the action of $Q$
on $H_*(\B N;k)$ is trivial, i.e. every element of $Q$ acts as the identity. 
Since $Q$ is abelian, the results of Dwyer and Greenlees 
\cite[propositions 5.3 and 6.9]{Dwyer Greenlees} show that $\chains_*(\B N;k)$ is a $k$-cellular
$kQ$-module. So, by lemma~\ref{pro: cellularization over a group-ring} for every $kQ$-module $M$, the
map $\cell^Q M \to M$ is a $k$-cellular approximation of $M$ as a $kG$-module.
\end{example}

\begin{example}
\label{exa: S3 doesnt have independence of base}
Here is an example where there is no independence of base property.
Let $\Sigma_3$ be the symmetric group on three elements. There is a short exact sequence of groups:
\[ C_3 \to \Sigma_3 \to C_2 \]
where $C_n$ is the cyclic group of $n$ elements.
Let $R$ be the group ring $\Int \Sigma_3$, let $S$ be the group-ring $\Int C_2$ and let $k$ be the $S$-module 
$\Int$ with trivial $C_2$-action. As noted in the previous example, $S \otimes_R k$ is equivalent to
the chains of the 
classifying space of $C_3$, namely $\chains_*(\B C_3 ; \Int)$. Recall $H_1 (\B C_3;\Int)\cong \Int/3$ 
(the group $H_1 (\B C_3;\Int)$ is $\pi_1 (\chains_*(\B C_3 ; \Int))$ in our usual notation).
It is easy to see that the $C_2$ action on $H_1 (\B C_3; \Int)$ is simply exchanging the two generators of $\Int/3$.

Let $I$ be the augmentation ideal of $S$, i.e. the kernel of the map $\Int C_2 \to \Int$.
Had $S \otimes_R k$ been $k$-cellular as an $S$-module, then $H_1 (\B C_3; \Int)$ would have been a nilpotent 
$S$-module by~\cite[proposition 6.11]{Dwyer Greenlees}. 
But it is easy to see that $H_1 (\B C_3; \Int)$ is not a nilpotent $S$-module.
Moreover, this argument also shows that $S \otimes_R k$ is not $k$-cellular even as an $R$-module.

\end{example}

We turn to discuss the applications of lemma~\ref{lem: Independence of Base} to topological groups.
We start with the following example.
\begin{example}
\label{exa:For a finite nilpotent space tildeB is a k-cellular pi_1(B)-space part 2}
Let $B$ be a finite connected nilpotent space and $k$ a commutative ring. The map $\kangr B \to \pi_1(B)$ of
topological groups induces a map $k[\kangr B] \to k[\pi_1(B)]$, we show this map has the strong independence of
base property.

Let $\tilde{B}$ be the universal cover of $B$. The following is a homotopy fibration sequence:
\[ \tilde{B} \to B \to K(\pi_1(B),1)\]
The map of $\sphere$-algebras we have in mind is the induced map of group rings: $k[\kangr B] \to k[\pi_1(B)]$.
As noted in example~\ref{exa:For a finite nilpotent space tildeB is a k-cellular pi_1(B)-space},
$\chains_*(\tilde{B};k)$ is a $k$-cellular $k[\pi_1(B)]$-module. Note that $k[\pi_1(B)]\otimes_{k[\kangr B]} k$
is equivalent to $\chains_*(\tilde{B};k)$. This implies strong independence of base 
(lemma~\ref{lem: Independence of Base}), i.e. for any $k[\pi_1(B)]$-module
$X$ the map $\cell^{\pi_1(B)} X \to X$ is a $k$-cellular approximation of $X$ as a $k[\kangr B]$-module
(strong independence of base).
The importance of this property will become clear in section~\ref{sec:Fibrations}.

\end{example}

The general property alluded to in 
example~\ref{exa:For a finite nilpotent space tildeB is a k-cellular pi_1(B)-space part 2} is given in the next
proposition. Note that if $F \to E \to B$ is a homotopy fibration sequence of spaces with $B$ a connected space,
then $F$ is equivalent to a $\kangr B$-space (we recall this standard construction
at the beginning of section~\ref{sec:Fibrations}). Hence we assume that $\chains_*(F;k)$ is a $k[\kangr B]$-space. 
\begin{proposition}
\label{pro: Strong independence for a fibration}
Let $F \to E \to B$ be a homotopy fibration sequence of spaces with $B$ and $E$ connected and let $k$ be a 
commutative $\sphere$-algebra. Suppose $\chains_*(F;k)$ is a $k$-cellular $k[\kangr B]$-module. Then
for every $k[\kangr B]$-module $X$ the map $\cell^{k[\kangr B]} X \to X$ is a $k$-cellular approximation
of $X$ as a $k[\kangr E]$-module.
\end{proposition}
\begin{proof}
From lemma~\ref{lem: Independence of Base} we see
it is enough to show that $k[\kangr B] \otimes_{k[\kangr E]} k$ is a $k$-cellular $k[\kangr B]$-module.
Hence it would suffice to show that $k[\kangr B] \otimes_{k[\kangr E]} k$ is equivalent to
$\chains_*(F;k)$ as $k[\kangr B]$-modules.

Consider the homotopy fibration sequence $\kangr B \to F \to E$. The Borel construction gives
an equivalence $ \kangr B \times_{\kangr E} \W \kangr E \simeq F$. Moreover, the map giving this equivalence is map 
of $\kangr B$-spaces. Thus we have an equivalence of $k[\kangr B]$-modules:
\[ \chains_*(\kangr B \times_{\kangr E} \W \kangr E;k) \simeq \chains_*(F;k)\]
To finish the proof we use the results of Elmendorf, Kriz, Mandell and May 
\cite[proposition IV.7.5 \& theorem IV.7.8]{EKMM}, which show that:
\[ \chains_*(\kangr B \times_{\kangr E} \W \kangr E;k) \simeq \chains_*(\kangr B;k) \otimes_{k[\kangr E]} 
\chains_*( \W \kangr E;k) \]
Since $\chains_*( \W \kangr E;k) \simeq k$, we are done.
\elaboration{ 
The Borel correspondence comes from the following fact. Consider the pullback diagram:
\[\xymatrixcompile{{\kangr X} \ar[d] \\
{\W \kangr X \times_X Y} \ar[r] \ar[d]& Y \ar[d] \\ {\W\kangr X} \ar[r] & X }\]
Thus $ \kangr X \to \W \kangr X \times_X Y \to Y$ is a principal fibration and hence
$W \kangr X \times_{\kangr X}(\W \kangr X \times_X Y) \simeq Y$. In fact this map is an isomorphism, as we will now
explain.

Let $G$ be a topological group, $A$ a right $G$-space and $B$ a left $G$ space. Suppose $B \to C$ is a
principal $G$-fibration and let $D \to C$ be a map. Then there is an isomorphism of simplicial
sets $A \times_G (B \times_C D) \cong (A \times_G B) \times_C D$. Thus we see that:
\[ (W \kangr X \times_{\kangr X}\W \kangr X) \times_X Y \cong X \times_X Y \cong Y\]

Now consider the diagram:
\[\xymatrixcompile{{\kangr E} \ar[d] \\
{X} \ar[r] \ar[d]& {\W\kangr E} \ar[d] \\ {W \kangr B \times_B E} \ar[r] & E }\]
We consider $\W \kangr E$ as a right $\kangr E$-space and $\W \kangr B$ as a left $\kangr B$-space.
Now $X=(\W \kangr B \times_B E) \times_E \W \kangr E \cong \W \kangr B \times_B\W \kangr E$.
Now we compute $X \times_{\kangr E} \W \kangr E$ using the isomorphism shown above:
\[ \cong \W \kangr B \times_B\W \kangr E \times_{\kangr E} \W \kangr E \cong
\W \kangr B \times_B E\]
This is clearly an equivalence of left $\kangr B$-spaces.
}
\end{proof}

\begin{example}
\label{exa: G to SU(N) to SU(N)/G part 1}
Let $G$ be a discrete finite group. Take an embedding of $G$ into the group $SU(n)$, there is always such an 
embedding for a large enough $n$. This embedding makes $SU(n)$ into a $G$-space.
Fix a commutative ring $k$.
We show that the principal fibration sequence 
\[ SU(n) \to SU(n)/G \to \B G\]
and the commutative ring $k$ satisfy the conditions of proposition~\ref{pro: Strong independence for a fibration}.
All we need is to show $\chains_*(\B\kangr SU(n);k)$ is a $k$-cellular $G$-module.
It is easy to see that the $G$-action on $\B \kangr SU(n) \simeq SU(n)$ comes from 
the embedding of groups $G \to SU(n)$. In~\cite[example 5.9]{Dwyer Greenlees Iyengar}, Dwyer Greenlees and
Iyengar show that this map makes $\chains_*(SU(n);k)$ into a $k$-cellular $kG$-module.
\end{example}

We end this section with an example of a topological group $G$ where the map $kG \to k[\pi_0(G)]$ does not
have strong independence of base property.
\begin{example}
Fix some natural number $n$ and consider the following fibration of pointed spaces:
\[ S^n \to \mathbb{RP}^n \to \B C_2\]
As in example~\ref{exa:For a finite nilpotent space tildeB is a k-cellular pi_1(B)-space part 2} we
set $R=\Int [\kangr \mathbb{RP}^n]$ (recall this is $\chains_*(\kangr \mathbb{RP}^n; \Int)$) and
$S = \Int [\kangr \B C_2]$. The weak equivalence of topological groups 
$\kangr \B C_2 \simeq C_2$ induces an equivalence between the derived category of $S$-modules and the derived
category of $\Int [C_2]$-modules (see~\cite[theorem III.4.2]{EKMM}). So we will work in the category of 
$\Int [C_2]$-modules instead of the category of $S$-modules. 

The action of $C_2$ on the $n$-sphere $S^n$ is by the antipodal map.
If $n$ is odd this map is homotopy equivalent to the identity. This implies that we are in the same situation as in 
example~\ref{exa:For a finite nilpotent space tildeB is a k-cellular pi_1(B)-space part 2}.
However, if $n$ is even, then the antipodal map induces a non-trivial action on the reduced homology
groups of the $n$-sphere $S^n$. Thus $\pi_n(\chains_*(S^n;\Int))$ is a non-trivial $C_2$-module,
the $C_2$-action being multiplication by $-1$. 
We will show this $C_2$-module, which we denote by $\tilde{\Int}$, is not $\Int$-cellular as a $\Int[C_2]$-module
(we consider $\Int$ to be a $\Int[C_2]$-module with the trivial action of $C_2$).

The ring $\Int [C_2]$ is commutative and its augmentation ideal is finitely generated,
so by~\cite[proposition 6.11]{Dwyer Greenlees}
a finitely generated $\Int [C_2]$-chain complex is $\Int$-cellular if and only if all of its homology groups are
nilpotent discrete $\Int [C_2]$-modules. It is easy to see that $\tilde{\Int}$ is not a nilpotent module and 
therefore $\chains_*(S^n;\Int)$ described above is not $\Int$-cellular as a $\Int [C_2]$-module.
\end{example} 

\section{Cellular Approximation in Nilpotent Groups}
\label{sec:Cellular Approximation in Nilpotent Groups}
Fix a prime number $p$ and let $k$ be a commutative $\sphere$-algebra such that $\pi_0(k)$ is an $\field_p$-algebra.
In this section we consider $k$-cellular approximation over a group-ring $kG$, where $G$ is a finite group.
When $G$ is a $p$-group we show that every $kP$-module is $k$-cellular, this is done in 
corollary~\ref{cor:Finite p-group}. 
For a finite nilpotent group $G$ we give a nice formula for $k$-cellular approximation over $kG$
in proposition~\ref{pro:Cellularization over a nilpotent group}. We further show, in 
proposition~\ref{pro: k proxy-small over kN}, that $k$ is proxy-small as a $kG$-module.

\subsection{Finite $p$-groups}
Recall that $\pi_0(k)$ is an $\field_p$-algebra.
This implies the multiplication by $p$ map: $k \xrightarrow{\cdot p} k$ is equal to the zero map
in $\Derived_k$.
Using this fact we will show that if $P$ is a finite $p$ group, then $k$ finitely builds $kP$. In particular 
we get that $k$ is a proxy-small $kP$-module and hence every $kP$-module is $k$-cellular. 
These results are obvious when $k$ is a commutative $\field_p$-algebra and the generalization to 
the case where $k$ is an $\sphere$-algebra is simple.

\begin{proposition}
\label{pro: k finitely builds kP}
Let $P$ be a finite $p$-group and $k$ a commutative $\sphere$-algebra such that $\pi_0(k)$ is an $\field_p$-algebra.
Then $k$ finitely builds $kP$ in the category of $kP$-modules and therefore $k$ is a proxy-small $kP$-module.
\end{proposition}

We start by proving the following lemma.
\begin{lemma}
\label{lem: triangle for a finite cyclic group}
Let $C$ be a finite cyclic group with a generator $g$ and let $k$ be any commutative $\sphere$-algebra, then there is a triangle: \[ \Sigma k \to K \to k\]
Where $K$ is the homotopy cofiber of the map $kC \xrightarrow{1-g} kC$.
\end{lemma}
Note that the construction of such a triangle is obvious when $k$ is a commutative ring. Indeed there is a
short exact sequence: \[ 0 \to k \to kC \xrightarrow{1-g} kC \to k \to 0\]
\begin{proof}[Proof of lemma~\ref{lem: triangle for a finite cyclic group}]
Let us begin with the case where $k$ is the sphere spectrum $\sphere$.
By \cite[IV.3.1]{EKMM} there is a map of $\sphere$-algebras $\sphere C \to H\pi_0(\sphere C)$.
Note that $H\pi_0(\sphere C)\cong H \Int C$.
Clearly ${H\Int C \otimes_{\sphere C} K}$ is equivalent to the cone of the map $\Int C \xrightarrow{1-g} \Int C$ 
and hence:
\[ \pi_i(H\Int C \otimes_{\sphere C} K)= \left\{ \begin{array}{cc}
  \Int & i=0,1 \\
  0 & \text{otherwise} \\
\end{array} \right . \]

The augmentation map $\sphere C \to \sphere$ induces a map $f:K \to \sphere$.
There is also a map $\sphere \to \sphere C$ of $\sphere C$-modules, which is defined by 
$1 \mapsto \sum_{x \in C} x$. One can easily check this is indeed a map of $\sphere C$-modules and it induces a map
$g:\Sigma \sphere \to K$.
\elaboration{ 
One can define this map in the following way: this is clearly a map is the category of $\sphere$-modules
(one considers $\sphere C$ simply as $\oplus_{c\in C} \sphere$).
Now, consider the category of $\sphere C$-modules simply as the category of $C$-diagrams over the category
of $\sphere$-modules. (Note that if one works with the $\sphere$-algebras as in \cite{EKMM} then one can simply use
$C$ as a discrete topological space, and thus $\sphere C$ is indeed the category of $C$-diagrams. However,
if one uses symmetric spectra, one must replace $C$ by the corresponding topological group, however $\sphere C$
would still be $\bigsqcup_C \sphere$ as an $\sphere$-module.) Then this map is indeed a map of $C$-diagrams, after 
giving $\sphere$ and $\sphere C$ the correct $C$-diagram structure.
}
Let $F$ be the homotopy fiber of the map $f$, then clearly:
\[ \pi_i(H\Int C \otimes_{\sphere C} F)= \left\{ \begin{array}{cc}
  \Int & i=1 \\
  0 & \text{otherwise} \\
\end{array} \right . \]
Since the composition $\sphere \to \sphere C \xrightarrow{1-g} \sphere C$ is zero
, one sees that the map $g$ induces a map $g':\Sigma \sphere \to F$. It is easy to check that the induced map
\[H\Int C \otimes_{\sphere C} \Sigma \sphere  \to H\Int C \otimes_{\sphere C} F\] is an equivalence.
Consider the mapping cone $\cone(g')$. This $\sphere C$-module is bounded below, simply because of it's
construction. 
In addition, $H\Int C \otimes_{\sphere C} \cone (g') \simeq 0$.
Hence, by \cite[IV.3.6]{EKMM}, $\cone(g') \simeq 0$ and the map $g'$ is an equivalence, yielding
the desired triangle:
\[ \Sigma \sphere \to K \to \sphere \]

Now suppose $k$ is any commutative $\sphere$-algebra. Note that: $k \simeq k \otimes_\sphere \sphere$
and $kC \simeq k \otimes_\sphere \sphere C$. Hence applying the functor $k \otimes_\sphere -$ to the triangle
$\Sigma \sphere \to K \to \sphere $ gives the triangle: $\Sigma k \to K \to k$.
\end{proof}

\begin{proof}[Proof of proposition~\ref{pro: k finitely builds kP}] 
We start with the case where $P$ is the cyclic group of $p$ elements.
Let $g$ be a generator of $P$ and let $K$ be the homotopy cofiber of the map $kP \xrightarrow{1-g} kP$.
By the previous lemma, $k$ finitely builds $K$. Denote by $K_n$ the homotopy cofiber of the map $kP 
\xrightarrow{(1-g)^n} kP$. Consider the following commutative diagram whose rows and columns are triangles:
\[ \xymatrixcompile{
{K_n} \ar[r] \ar[d] & {K_{n+1}} \ar[r] \ar[d] & {K} \ar[d] \\
{kP} \ar@{=}[r] \ar[d]^{(1-g)^n} & {kP} \ar[r] \ar[d]^{(1-g)^{n+1}} & 0 \ar[d] \\
{kP} \ar[r]^{1-g} & {kP} \ar[r] & {\Sigma K} }\]
From the top row we see, by induction on $n$, 
that $K_n$ is finitely built from $K$. In particular $K_p$ is finitely built by $K$.
Since the map $(1-g)^p$ is equivalent to
the zero map, the module $K_p$ is equivalent to $kP\oplus \Sigma^{-1} kP$. Therefore $K_p$ finitely builds
$kP$ and hence $k$ finitely builds $kP$. This completes the case where $P$ is the cyclic group with $p$-elements.

Now suppose $P$ is any finite $p$-group. We prove the proposition by induction on the order of $P$.
Since $P$ is a finite $p$-group, there is a short exact sequence of groups:
\[ N \hookrightarrow P\twoheadrightarrow C_p \]
Where $C_p$ is the cyclic group with $p$-elements.
By the induction assumption, $k$ finitely builds $kN$ in the category of $kN$-modules.
By applying the functor $kP \otimes_{kN} -$ to the recipe of constructing $kN$ from $k$, we see that
$kP \otimes_{kN} k \simeq kC_p$ finitely builds $kP \otimes_{kN} kN \simeq kP$ in the category of $kP$-modules.

We saw that $k$ finitely builds $kC_p$ in the category of $kC_p$-modules. Hence 
$k$ also finitely builds $kC_p$ in the category of $kP$-modules. We conclude $k$ finitely builds $kP$ 
in the category of $kP$-modules.

Finally, $kP$ is a Koszul-complex for $k$, since $kP$ is a small $kP$-module, $k$ finitely builds $kP$, and
clearly $kP$ builds $k$.
\end{proof}

Here is the obvious corollary of proposition~\ref{pro: k finitely builds kP}.
\begin{corollary}
\label{cor:Finite p-group}
Let $P$ be a finite $p$-group and $k$ a commutative $\sphere$-algebra such that $\pi_0(k)$ is an $\field_p$-algebra.
Then every $kP$-module is $k$-cellular.
\end{corollary}

\subsection{Finite Nilpotent Groups}
Let $N$ be a finite nilpotent group. As before, $p$ is a fixed prime number and $k$ is a commutative 
$\sphere$-algebra such that $\pi_0(k)$ is an $\field_p$-algebra.
Let $P$ be a $p$-Sylow subgroup of $N$ and let $H=N/P$. Since $N$ is nilpotent, it is isomorphic to the
group $P \times H$.

\begin{proposition}
\label{pro:Cellularization over a nilpotent group}
Let $N$, $P$, $H$ and $k$ be as above. Then for any $N$-module $X$, the map:
$\Hom_{kN}(kP,X)\to X$ (induced by the map $kN \to kP$), is the $k$-cellular approximation of
$X$ as an $N$-module.
\end{proposition}
\begin{proof}
There is a short exact sequence of groups:
$1\to H \to N \to P \to 1$.
As we saw in proposition~\ref{pro: k finitely builds kP}, $k$ builds $kP$. Hence $k$-cellular approximations
are the same as $kP$-cellular approximations.

Since $P$ is a $p$-Sylow subgroup, the order of $H$ is prime to $p$. This implies the map:
$k \xrightarrow{\cdot |H|} k$ (multiplication by the order of $H$) is an equivalence. Let $b:k \to kH$ be the map
defined (up to homotopy) by sending $1$ to $\Sigma_{h\in H} h$.
The composition $k \to kH \to k$ is the multiplication by $|H|$ map, hence an equivalence.
Thus the augmentation map $kH \to k$ has a right inverse, showing $k$ is a retract of $kH$ 
(in the derived category $\Derived_{kH}$). Applying the functor $kN \otimes_{kH} -$ to these maps yields two maps:
$kP \to KN \to kP$ whose composition is the identity (in $\Derived_{kN}$).
Note that the right map is the one induced by the
map $N \to P$ of groups. Now we use lemma~\ref{lem: k is a retract of R}, showing:
$\Hom_{kN} (kP,X) \to X$ is a $k$-cellular approximation of $X$ as a $kH$-module.
\end{proof}
\begin{remark}
\label{rem: strict fixed points}
Note that the map $\Hom_{kN}(kP,X) \to X$ is equivalent to the map $\Hom_{kH}(k,X)\to X$ by well known adjunctions.
Suppose that $k = \field_p$, then $k$ is a projective $H$-module (because it is a retract of $kH$).
This implies the derived functor: $\Hom_{kH}(k,-)$ is equivalent to the non-derived version.
So, for a $kN$-chain complex $X$, $\cell^N X$ is the fixed points of the $H$-action on $X$, 
commonly denoted $X^H$.
\end{remark}

In fact, the proof of proposition~\ref{pro:Cellularization over a nilpotent group} implies a slightly more
general result:
\begin{proposition}
Let $H \to G \to P$ be a short exact sequence of finite groups, with $P$ being a $p$-group.
Suppose that the order of $H$ is prime to $p$ and that $\pi_0(k)$ is an $\field_p$-algebra, then
for any $kG$-module $X$ the natural map: $\Hom_{kG}(kP,X) \to X$ is a $k$-cellular approximation of $X$ as
a $kG$-module.
\end{proposition}

We also show that $k$ is proxy-small as a $kN$-module. This will later enable us to use Dwyer, Greenlees and
Iyengar's formula for cellular approximation \cite[theorem 4.10]{Dwyer Greenlees Iyengar}.
\begin{proposition}
\label{pro: k proxy-small over kN}
Let $N$ be a finite nilpotent group and $k$ a commutative $\sphere$-algebra such that $\pi_0(k)$ is an 
$\field_p$-algebra. Then $k$ is a proxy-small $kN$-module.
\end{proposition}
\begin{proof}
As before let $P$ be the $p$-Sylow subgroup of $N$ and let $H=N/P$ so $N \cong P \times H$.
We will show that $kP$ is a Koszul complex for $k$ over $kN$.
From the proof of proposition~\ref{pro:Cellularization over a nilpotent group} above we see that $k$ is a 
retract of $kH$ (in $\Derived_{kN}$). This implies $kN \otimes_{kH} k \simeq kP$ is a retract
of $kN \otimes_{kH} kH \simeq kN$. Hence $kP$ is small as an $N$-module.

From proposition~\ref{pro: k finitely builds kP} we see $k$ finitely builds $kP$ over $kP$ and therefore $k$ also
finitely builds $kP$ over $kN$. Clearly $kP$ builds $k$ over $kN$, since $k$ is a $kP$-module.
Hence $kP$ is a Koszul-complex for $k$ over $kN$.
\end{proof}

\section{Relation with the target of the Eilenberg-Moore Spectral Sequence}
\label{sec:Fibrations}


Fix a fibration sequence $F \to E \to B$, where $E$ and $B$ are always assumed to be
connected spaces. Let $\W \kangr B$ be a contractible space with a free $\kangr B$-action 
(such as the one described by Rothenberg and Steenrod in \cite{RothenbergSteenrod}).
The pullback of the diagram:
\[\xymatrixcompile{{\Box} \ar@{-->}[r] \ar@{-->}[d]& E \ar[d] \\ {\W\kangr B} \ar[r] & B }\]
is a (right) $\kangr B$-space equivalent to $F$. So, without loss of generality, we will assume $F$ is a $\kangr B$-space. Let $k$ be a commutative $\sphere$-algebra. In this section the chains of spaces are always taken  
with coefficients in $k$, therefore we omit $k$ from the notation $\chains_*(-;k)$.
We use $R$ to denote the group-ring $k[\kangr B]$, i.e. $\chains_*(\kangr B;k)$.
So $\chains_*(F)$ and $\chains^*(F)$ are both $R$-modules. 


The results of Dwyer, Greenlees and Iyengar imply that, under some conditions on the space $B$,
the natural map $\chains^*(E) \otimes_{\chains^*(B)} k \to \chains^*(F)$ is a $k$-cellular approximation of
$\chains^*(F)$ over $R$ (see~\cite[theorem 4.10]{Dwyer Greenlees Iyengar}). These conditions
on $B$ are specified in lemma~\ref{lem: fibration EMSS gives cellularization} below. Applying the machinery we have
set up in the previous section to the $k$-cellular approximation of $\chains^*(F)$ over $R$, we obtain the following 
theorems.

\begin{theorem}
\label{the: Convergence of EMSS over K(P,1)}
Fix a prime $p$. Let $N$ be a finite nilpotent group and let $P\subseteq N$ be the $p$-Sylow subgroup of $N$,
so $N \cong P \times H$ with the order of $H$ being prime to $p$. 
Let $F \to E \to \B N$ be a homotopy fibration sequence over the classifying space of $N$,
with $E$ being a connected space.
If $k$ is any commutative $\sphere$-algebra such that $\pi_0(k)$ is an $\field_p$-algebra, then:
\[ \chains^*(E;k) \otimes_{\chains^*(B;k)} k \simeq \chains^*(F_{h(H)})\]
Where $F_{h(H)}$ is the homotopy orbit space of $F$ with respect to the $H$-action on $F$.
In particular, if $N=P$ then \[ \chains^*(E;k) \otimes_{\chains^*(B;k)} k \simeq \chains^*(F;k) \]
\end{theorem}

\begin{theorem}
\label{the: Convergence of EMSS to fixed points over finite nilpotent space}
Fix a prime $p$. Let $B$ be a finite connected nilpotent space with a finite fundamental group $N=\pi_1(B)$.
Let $P\subseteq N$ be the $p$-Sylow subgroup of $N$, so $N \cong P \times H$.
Let $F \to E \to B$ be a homotopy fibration sequence over $B$, where $E$ is a connected space. Then:
\[ \tor_{-n}^{\chains^*(B;\field_p)}(\chains^*(E;\field_p), \field_p) = H^n(F;\field_p)^{H}\]
Where $H^n(F;\field_p)^{H}$ are the fixed points of the $H$-action on $H^n(F;\field_p)$.
\end{theorem}

\begin{remark}
The tensor product $\chains^*(E) \otimes_{\chains^*(B)} k$ needs, perhaps, some clarification.
The maps $E \to B$ and $pt \to B$ induce maps of $k$-algebras $\chains^*(B) \to \chains^*(E)$ and
$\chains^*(B) \to k$. These maps make $\chains^*(E)$ and $k$ into $\chains^*(B)$-bimodules. Considering
$\chains^*(E)$ as a right $\chains^*(B)$-module and $k$ as a left $\chains^*(B)$-module, we form
the tensor product: $\chains^*(E) \otimes_{\chains^*(B)} k$.
\end{remark}

\begin{remark}
The Eilenberg-Moore spectral sequence is of the form:
\[ E^2_{p,q}=\tor_{p,q}^{\pi_*(\chains^*(B))}(\pi_*(\chains^*(E)),\pi_*(k)) \ \Rightarrow \ 
\pi_{p+q} (\chains^*(E) \otimes_{\chains^*(B)} k)\]
The version of the Eilenberg-Moore spectral sequence we consider is the one given in \cite[IV.4.1]{EKMM}.
Its convergence properties are also described there and elsewhere in the literature and will not
be treated here. However, note that theorems \ref{the: Convergence of EMSS over K(P,1)} and 
\ref{the: Convergence of EMSS to fixed points over finite nilpotent space}
do describe the $E^\infty$ term of this spectral sequence.  
\end{remark}

As noted in the introduction, the first result (theorem~\ref{the: Convergence of EMSS over K(P,1)})
generalizes a result of Kriz from~\cite{Kriz}, showing convergence of the Eilenberg-Moore 
cohomology spectral sequence for a fibration over $\B (\Int/p)$, with coefficients in Morava $K$-theory.
The second result (theorem~\ref{the: Convergence of EMSS to fixed points over finite nilpotent space})
is the dual to a result of Dwyer \cite{Dwyer} concerning convergence of the
Eilenberg-Moore homology spectral sequence and in fact follows from it.

We start by recalling some known properties of fibrations and $\sphere$-algebras.
\begin{lemma}
\label{lem: properties of fibrations 1}
Let $F \to E \to B$ be a fibration as above and let $R=k[\kangr B]$, then:
\begin{enumerate}
\item
$\Hom_R(k,k) \simeq \chains^*(B)$ as $\sphere$-algebras.
\item
$\chains_*(F) \otimes_R k \simeq \chains_*(E)$.
\item
The modules $\Hom_R(k,\chains^*(F))$ and $\chains^*(E)$ are equivalent as right $\chains^*(B)$-modules.
\end{enumerate}
\end{lemma}
\begin{proof}
First, recall that $\W \kangr B$ is a contractible free $\kangr B$-space. The term "free" is as defined by
Rothenberg and Steenrod in \cite{RothenbergSteenrod}. It is easy to see that 
$\chains_*(\W \kangr B ;k)$ is a cell $R$-module (see \cite[III.2.1]{EKMM}) that is equivalent to $k$ as $R$-modules.

For the first equivalence, note that $\W \kangr B \times_{\kangr B} pt$ is equivalent to the classifying space of $\kangr B$ which is $B$. Therefore, by~\cite[proposition IV.7.5 \& theorem IV.7.8]{EKMM},
we see $k \otimes_R k \simeq \chains_*(B;k)$. This equivalence implies:
\begin{align*}
\chains^*(B) & \simeq \Hom_k(\chains_*(B;k),k) \simeq \Hom_k (k \otimes_R k,k) \\
& \simeq \Hom_R(k,\Hom_k(k,k)) \simeq \Hom_R(k,k)
\end{align*}
This equivalence is indeed an equivalence of $\sphere$-algebras. This was noted by
Dwyer and Wilkerson in~\cite[proof of lemma 2.10]{Dwyer Wilkerson} and also by Dwyer, Greenlees and Iyengar
in \cite[4.22]{Dwyer Greenlees Iyengar}.
We sketch the argument here. Consider the following maps of $\sphere$-algebras, all of which are equivalences:
\begin{align*}
\chains^*(B)&= \Hom_k(\chains_*(B;k),k) \xrightarrow{\simeq} 
\Hom_k(\chains_*(\W \kangr B \times_{\kangr B} pt);k)\\ &\cong
\Hom_k(\chains_*(\W \kangr B) \otimes_R k ;k) \cong
\Hom_R(\chains_*(\W \kangr B) ,\chains_*(\W \kangr B))
\end{align*}
The first isomorphism follows from the structure of $\W \kangr B$ as a free $\kangr B$-space. The second 
isomorphism is \cite[proposition III.6.3]{EKMM}.

For the second equivalence, recall the homotopy orbit space of the $\kangr B$-space $F$ is
$F_{h \kangr B} = F \times_{\kangr B} \W \kangr B$. There is a well known equivalence $E \simeq F_{h \kangr B}$
(this is sometimes called the Borel correspondence).
Using again the results of~\cite[proposition IV.7.5 \& theorem IV.7.8]{EKMM}, we see that 
$ \chains_*(F) \otimes_R k \simeq \chains_*(E)$ (we consider $F$ as a right $\kangr B$-space, thus
$\chains^*(F)$ becomes a left $R$-module).

This last equivalence implies that:
\begin{align*}
\Hom_R(k,\chains^*(F)) & \simeq \Hom_R(k,\Hom_k(\chains_*(F),k)) \\ 
&\simeq \Hom_k(\chains_*(F)\otimes_R k,k) \simeq \chains^*(E)
\end{align*}
Both $\Hom_R(k,\chains^*(F))$ and $\chains^*(E)$ are right $\chains^*(B)$-modules.
The module $\Hom_R(k,\chains^*(F))$ is naturally a right $\Hom_R(k,k)$-module by the composition pairing:
\[\Hom_R(k,\chains^*(F)) \otimes_k \Hom_R(k,k) \to \Hom_R(k,\chains^*(F))\]
(see \cite[III.6.12]{EKMM}). The map $\chains^*(B) \xrightarrow{\simeq} \Hom_R(k,k)$ makes
$\Hom_R(k,\chains^*(F))$ into a right $\chains^*(B)$-module.
For $\chains^*(E)$, it is a right $\chains^*(B)$-module by the map of $\sphere$-algebras
$\chains^*(B) \to \chains^*(E)$.

The equivalence $\Hom_R(k,\chains^*(F))\simeq \chains^*(E)$ is an equivalence of right $\chains^*(B)$-modules.
An outline of proof of this fact, for the case $k=H \field_p$,
can be found in~\cite{Dwyer Wilkerson} - in the course of the proof of lemma 2.10 there.
These arguments of Dwyer and Wilkerson from~\cite{Dwyer Wilkerson} readily generalize to an arbitrary
$\sphere$-algebra $k$, showing the modules $\Hom_R(k,\chains^*(F))$ and $\chains^*(E)$ are equivalent as right 
$\chains^*(B)$-modules.
\end{proof}

As mentioned at the beginning of this section,
results of Dwyer, Greenlees and Iyengar from~\cite{Dwyer Greenlees Iyengar} show, 
under some conditions, that the natural map
$\chains^*(E) \otimes_{\chains^*(B)} k \to \chains^*(F)$ is a $k$-cellular approximation of $\chains^*(F)$ over $R$.
These conditions are summarized in the following lemma.
\begin{lemma}
\label{lem: fibration EMSS gives cellularization}
If any one of the three conditions below holds, then $k$ is proxy-small as an $R$-module and
the map $\chains^*(E) \otimes_{\chains^*(B)} k \to\chains^*(F)$ (defined in $\Derived_R$) is a $k$-cellular 
approximation of $\chains^*(F)$ as an $R$-module.
\begin{enumerate}
\item
$k$ is a commutative $\sphere$-algebra and $B$ is a finite space.
\item  
$k$ is a commutative $\sphere$-algebra such that $\pi_0(k)$ is an $\field_p$ algebra
and $B$ is equivalent to an Eilenberg-Mac Lane space of type $K(N,1)$ where $N$ is a finite nilpotent group.
\item
$k$ is a commutative ring (in the classical sense) and $B\simeq K(G,1)$, where $G$ is a finite group.
\end{enumerate}
\end{lemma}
\begin{proof}
We first explain the origin of the map $\chains^*(E) \otimes_{\chains^*(B)} k \to\chains^*(F)$.
By the previous lemma (\ref{lem: properties of fibrations 1}) there is an equivalence:
\[\chains^*(E) \otimes_{\chains^*(B)} k \simeq \Hom_R(k,\chains^*(F)) \otimes_{\End_R(k)} k\]
This is an equivalence of $R$-modules, since the $R$-module structure on each of the tensor products
is induced by the $R$-module structure of $k$.
There is a natural map $\Hom_R(k,\chains^*(F)) \otimes_{\End_R(k)} k \to \chains^*(F)$ of $R$-modules 
(adjoint to the identity map $\Hom_R(k,\chains^*(F)) \xrightarrow{id} \Hom_R(k,\chains^*(F))$ of
right $\End_R(k)$-modules). This is the map referred to in the statement of the lemma.

If $k$ is proxy-small as an $R$-module, then, by \cite[theorem 4.10]{Dwyer Greenlees Iyengar},
the map: $\Hom_R(k,\chains^*(F)) \to \chains^*(F)$ is a $k$-cellular approximation of $\chains^*(F)$ as
an $R$-module.
Every one of conditions above ensures $k$ is a proxy-small $R$-module. The first condition 
by~\cite[proposition 5.3]{Dwyer Greenlees Iyengar}, the second condition by
proposition~\ref{pro: k proxy-small over kN} and the third condition by~\cite[example 5.9]{Dwyer Greenlees Iyengar}.
\end{proof}

\begin{example}
\label{exa: G to SU(N) to SU(N)/G part 2}
We continue with example~\ref{exa: G to SU(N) to SU(N)/G part 1}. Recall $G$ is a finite group and
$G \to SU(n)$ is an embedding of groups. As in example~\ref{exa: G to SU(N) to SU(N)/G part 1} we take
$k$ to be a commutative ring. Consider the fibration sequence:
\[ G \to SU(n) \to SU(n)/G\]
The base space of this fibration is a finite space, thus satisfying the first condition of 
lemma~\ref{lem: fibration EMSS gives cellularization}. Therefore the map:
\[ \chains^*\big(SU(n)\big) \otimes_{\chains^*(SU(n)/G)} k \to kG\]
is a $k$-cellular approximation of $R$-modules, $R$ being the group ring $k\big[\kangr\big(SU(n)/G\big)\big]$.
In example~\ref{exa: G to SU(N) to SU(N)/G part 1} we
saw there is a strong independence of base (lemma~\ref{lem: Independence of Base})
with respect to the map $R \to kG$ and the module 
$k$\elaboration{ (this map is induced by the isomorphism $\pi_0(\kangr(SU(n)/G)) \cong G$)}.
In particular $\cell_k^R (kG)$ is equivalent to $\cell^G (kG)$. We conclude that:
\[  \chains^*\big(SU(n)\big) \otimes_{\chains^*(SU(n)/G)} k \simeq \cell^G (kG)\]
This implies that the target of the Eilenberg-Moore cohomology spectral sequence for the fibration
$G \to SU(n) \to SU(n)/G$ is $\cell^G kG$.
For example, if $G$ is an abelian group, then the results of Dwyer and Greenlees \cite{Dwyer Greenlees} imply
that
\[ \pi_{-q} \bigg(\chains^*\big(SU(n)\big) \otimes_{\chains^*(SU(n)/G)} k \bigg) \cong H_I^q(kG)\]
where $I$ is the augmentation ideal of $kG$ and $H_I^*(-)$ denotes the local cohomology groups
with respect to $I$.
\end{example}

Before continuing to the proofs of theorems~\ref{the: Convergence of EMSS over K(P,1)} 
and~\ref{the: Convergence of EMSS to fixed points over finite nilpotent space}, we demonstrate the use of
cellular approximations by giving a different proof to the following result of Dwyer from~\cite{DwyerStrong}.
\begin{proposition}
\label{pro: Different proof of Dwyer's strong convergence}
Let $k$ be a commutative ring and let $F \to E \to B$ be a fibration sequence where $E$ and $B$ are connected and $B$ is either a finite space or the classifying space of some finite group.
Suppose $\pi_1(B)$ acts nilpotently on the cohomology groups $H^n(F;k)$, i.e. for every $n \geq 0$ and 
$x \in H^n(F;k)$ there exists some $m \geq 0$ such that $(1-g)^m x =0$ for all $g \in \pi_1(B)$.
Then there is an equivalence of $k\kangr B$-modules:
\[ \chains^*(E) \otimes_{\chains^*(B)} k \simeq \chains^*(F)\]
\end{proposition}
\begin{proof}
As usual, we use $R$ to denote the $\sphere$-algebra $k[\kangr B]$.
From lemma~\ref{lem: fibration EMSS gives cellularization} we see it is sufficient to show that
$\chains^*(F)$ is a $k$-cellular $R$-module. Because $R$ is connective and
$\chains^*(F)$ is bounded above, $\chains^*(F)$ is built from it's homotopy groups
(see the proof of proposition~\ref{pro: coconnective module over a connective algebra A is A_0-cellular}).
In other words, $\chains^*(F)$ is built by the $k [\pi_1(B)]$-modules $\{ H^n(F;k) \}_{n \geq 0}$.
So, it is enough to show that for every $n$, $H^n(F;k)$ is a $k$-cellular $R$-module.

Fix some $n$. The $k[ \pi_1(B)]$-module $H^n(F;k)$ can be written as an increasing union 
of nilpotent $k[\pi_1(B)]$-modules: $H^n(F;k)=\bigcup_{i\geq1} N_i$ where $N_i$ is a nilpotent $k[\pi_1(B)]$-module
of class $i$. It is easy to see, using induction on $i$, that each $N_i$ is a $k$-cellular $k[\pi_1(B)]$-module.
The $k[\pi_1(B)]$-module $H^n(F;k)$ is equivalent to the homotopy colimit of the telescope
$N_1 \to N_2 \to \cdots$. Since the homotopy colimit of $k$-cellular modules
is $k$-cellular (see, for example~\cite{Farjoun}), we conclude $H^n(F;k)$ is $k$-cellular as a $k [\pi_1(B)]$-module.
But this implies $H^n(F;k)$ is $k$-cellular also as an $R$-module.
\end{proof}

\begin{proof}[Proof of theorem~\ref{the: Convergence of EMSS over K(P,1)}]
By \cite[theorem III.4.2]{EKMM}, the equivalence $\kangr B \simeq N$ induces an equivalence of the derived 
categories: $\Derived_R \simeq \Derived_{kN}$.
This means we can work in the category of $kN$-modules instead of $R$-modules.
To do that we replace $F$ by an equivalent $N$-space, which we will also denote 
$F$\elaboration{ (this is done by taking $N \times_{\kangr B} F$ instead of $F$)}.
Since $\chains^*(E) \otimes_{\chains^*(B)} k \simeq \cell^N \chains^*(F)$,
the result follows from proposition~\ref{pro:Cellularization over a nilpotent group} by noting that: 
\[\Hom_{k[H]}(k,\chains^*(F)) \simeq \Hom_k(\chains_*(F) \otimes_{k[H]} k,k)
\simeq \chains^*(F_{h(H)})\]
The last equivalence: $\chains_*(F) \otimes_{k[H]} k \simeq F_{h(H)}$, follows 
from two results of Elemendorf, Kriz, Mandell and May \cite[proposition IV.7.5 \& theorem IV.7.8]{EKMM}.
\end{proof}

Our next goal is to prove theorem~\ref{the: Convergence of EMSS to fixed points over finite nilpotent space}.
For that purpose, we will decompose the topological group $\kangr B$ into the following homotopy fibration 
sequence:
\[ \kangr \tilde{B} \to \kangr B \to \kangr K(\pi_1(B),1)\]
where $\tilde{B}$ is the universal cover of $B$ and $K(\pi_1(B),1)$ is the appropriate Eilenberg-Mac Lane space
(i.e. the classifying space of $\pi_1(B)$).
\elaboration{
This short exact sequence comes from applying the Kan loop-group functor to the fibration of connected
spaces: $ \tilde{B} \to B \to K(\pi_1(B),1)$.}

Let $Q$ denote the topological group $\kangr K(\pi_1(B),1)$ and set $G=\pi_1(B)$.
Instead of working with $kQ$-modules, we would much rather work with $kG$-modules.
The following observation shows that for all of our purposes we can replace $kQ$-modules by $kG$-modules.
The map $\mu:Q \to \pi_0(Q)=G$ is a weak equivalence and therefore the induced map $\mu:\chains_*(Q;k) \to kG$
of $\sphere$-algebras is also a weak equivalence.
By \cite[theorem III.4.2]{EKMM}, the functor $ \Derived_{kG} \xrightarrow{\mu^*} \Derived_{kQ}$, induced by $\mu$,
is an equivalence of these derived categories.
This implies we can work over $kG$ instead of $kQ$. Moreover, the map $\varphi:R \to kG$, induced by
the map $\kangr B \to \pi_0(\kangr B)=G$, is equal to the composition
$ R \xrightarrow{\psi} kQ \xrightarrow{\mu} kG$. Hence the functors: $\psi^* \mu^*$ and $\varphi^*$ are equal.

Here is the main idea of the proof of 
theorem~\ref{the: Convergence of EMSS to fixed points over finite nilpotent space}.
Recall $F \to E \to B$ is a fibration with $B$ a finite, connected, nilpotent space and $k$ is the field $\field_p$. 
From the proof of proposition~\ref{pro: coconnective module over a connective algebra A is A_0-cellular} we see that
$\chains^*(F)$ is built from the $k[\pi_1(B)]$-modules: $\{H^n(F;k)\}_{n \geq 0}$.
We will show that the $k$-cellular approximation of $\chains^*(F)$ (which is $\chains^*(E) \otimes_{\chains^*(B)} k$)
is built from the $k[\pi_1(B)]$-modules: 
$\{\cell^\pi_1(B) H^n(F;k)\}_{n\geq 0}$. Then we apply the results of 
section~\ref{sec:Cellular Approximation in Nilpotent Groups} to compute the modules 
$\cell^{\pi_1(B)} H^n(F;k)$.

\begin{proof}[Proof of theorem~\ref{the: Convergence of EMSS to fixed points over finite nilpotent space}]
Set $k=\field_p$, we remind the reader that
$R$ denotes the group ring $k[\kangr B]$, $N$ is the fundamental group $\pi_1(B)$, 
$P$ is a $p$-Sylow  subgroup of $N$ and $H=N/P$.
By lemma~\ref{lem: fibration EMSS gives cellularization}, 
$\chains^*(E) \otimes_{\chains^*(B)} k$ is a $k$-cellular approximation of $\chains^*(F)$ as an $R$-module.
So it is enough to show $\cell^R_k \chains^*(F)$ has the desired homotopy groups. 

Example~\ref{exa:For a finite nilpotent space tildeB is a k-cellular pi_1(B)-space}
shows that $\chains_*(\tilde{B})$ is a $k$-cellular $kN$-module.
So by the proposition~\ref{pro: Strong independence for a fibration},
for every $N$-module $X$, the map $\cell^N X \to X$ is a $k$-cellular approximation of $X$ as an $R$-module.
Proposition~\ref{pro:Cellularization over a nilpotent group} shows the map:
$\Hom_{kN}(kP,X) \to X$ is a $k$-cellular approximation of $X$ as a $kN$-module.
In fact, as noted in remark~\ref{rem: strict fixed points}, $\Hom_{kN}(kP,X)$ is the strict fixed points of the $H$
action on $X$, denoted $X^{H}$. 

We use now the notation of Dwyer and Greenlees from \cite{Dwyer Greenlees}, this is the same notation
used in the proof of proposition~\ref{pro: coconnective module over a connective algebra A is A_0-cellular}.
In this notation $H^n(F;k)\cong \chains^*(F)\langle -n,-n\rangle $ and
$ \cell_k^R \big(\chains^*(F)\langle -n,-n\rangle\big)  \simeq H^n(F;k)^{H}$.
We now prove, by induction on $i$, that for every $n$ such that $0\leq n \leq i$: 
\begin{equation}
\label{equ: the homotopy groups of cell F}
\pi_{-n}\cell_k^R \big(\chains^*(F)\langle -i,0\rangle \big) \cong H^n(F;k)^{H}
\end{equation}
For the induction step, consider the triangle:
\[ \chains^*(F)\langle -i,0\rangle  \to \chains^*(F)\langle -i-1,0\rangle \to \chains^*(F)\langle -i-1,-i-1\rangle \]
taking $k$-cellular approximation of this triangle yields the following triangle:
\[ \cell_k^R \big(\chains^*(F)\langle -i,0\rangle \big) \to \cell_k^R \big(\chains^*(F)\langle -i-1,0\rangle \big)
 \to \cell_k^R\big(\chains^*(F)\langle -i-1,-i-1\rangle \big) \]
The long exact sequence of homotopy groups, for this last triangle, shows that
\eqref{equ: the homotopy groups of cell F} holds for $0 \leq n \leq i+1$, thus completing the induction.

From the proof of proposition~\ref{pro: coconnective module over a connective algebra A is A_0-cellular} we see that:
$\chains^*(F) \simeq \hocolim_i \chains^*(F)\langle -i,0\rangle $. Next we show that:
\begin{equation}
\label{equ: cell commutes with hocolim}
    \cell_k^R \chains^*(F) \simeq \hocolim_i \cell_k^R\big(\chains^*(F)\langle -i,0\rangle \big)
\end{equation}
By a result of Dwyer and Greenlees \cite[proposition 4.3]{Dwyer Greenlees},
the $k$-cellular approximation of an $R$-module $X$ is given by the map:
$\cell_k^R (R) \otimes_R X \to X$, whenever $k$ is proxy-small. The functor $\cell_k^R(R) \otimes_R -$ is a left 
adjoint and therefore commutes with homotopy colimits. Hence taking $k$-cellular approximations commutes with 
homotopy colimits and the desired equivalence in \eqref{equ: cell commutes with hocolim} follows.

Finally, from \eqref{equ: the homotopy groups of cell F} and \eqref{equ: cell commutes with hocolim} above
it easy to see that:
\[\pi_{-n}\big(\cell_k^R \chains^*(F;k)\big) \cong 
\pi_{-n} \bigg(\cell_k^R\big(\chains^*(F;k)\langle -n,0\rangle \big)\bigg) \cong H^n(F;k)^{H}\]
\end{proof}

\elaboration{
\begin{remark}
The proof of theorem~\ref{the: Convergence of EMSS to fixed points over finite nilpotent space} works just as
well if we replace $\field_p$ with any commutative, connective and bounded above $\sphere$-algebra $k$ such that 
$\pi_0(k)$ is an $\field_p$-algebra. Because, when $k$ is connective we still have Postnikov sections and
if $k$ is bounded above then $\chains^*(F)$ is also bounded above. Note that $\chains_*(\tilde{B};k)$ is still
a $k$-cellular $N$-module (see example~\ref{exa:For a finite nilpotent space tildeB is a k-cellular pi_1(B)-space}).
The only difficulty is in showing that $k \otimes_{k[H]} H \pi_i(\chains^*(F;k)) \simeq 
H(\pi_i(\chains^*(F;k))^{H})$. This can be seen from the fact that $k$ is a retract of $k[H]$.
\end{remark}}

\section{A spectral sequence}
\label{sec: Applications}
As before, let $F \to E \to B$ be a fibration sequence.
The proof of theorem~\ref{the: Convergence of EMSS to fixed points over finite nilpotent space} shows
that when $B$ is a finite nilpotent space and $k$ is a commutative ring 
then $\cell_k^R\chains^*(F)$ is built from the modules $\cell^{k[\pi_1(B)]}_k H^n(F;k)$.
We use this observation to construct a spectral sequence converging (conditionally) to $\cell_k^R \chains^*(F;k)$, 
where $k$ is any commutative ring. Below, 
in~\ref{sub: Fibrations over a nilpotent space with cyclic fundamental group}, we demonstrate the use of this 
spectral sequence for the case where $\pi_1(B)$ is a cyclic group.
\begin{proposition}
\label{pro: Spectral Sequence}
Let $F \to E \to B$ be a fibration where $E$ and $B$ are connected and let $k$ be a commutative ring.
Suppose $B$ is a finite nilpotent space and set $N=\pi_1(B)$. There exists a spectral sequence:
\[ E^1_{p,q} = \pi_{2p+q}(\cell^N H^p(F;k)) \ \Rightarrow \ \pi_{p+q}( \cell_k^R \chains^*(F;k)) \]
where $R$ denotes the $\sphere$-algebra $k[\kangr B]$.
Convergence is as in the usual case for the spectral sequence of the homotopy colimit of spectra.
\end{proposition}
\begin{proof}
The $\sphere$-algebra $R=k[\kangr B]$ is connective and bounded-above.
Connectivity of $R$ implies the existence of Postnikov sections in the category of $R$-modules
(see~\cite[lemma 3.2]{Dwyer Greenlees Iyengar}). Since $k$ is bounded above, the $R$-module $\chains^*(F)$ is
also bounded above. So, the same arguments as in the proof of 
theorem~\ref{the: Convergence of EMSS to fixed points over finite nilpotent space} show that:
\[ \cell_k^R \chains^*(F) \simeq \hocolim_p \cell_k^R (\chains^*(F)\langle -p,0 \rangle)\]
This gives a spectral sequence whose $E^1$-term is:
\[ E^1_{p,q}=
\pi_{p+q} \bigg(\cone \Big(\cell^R_k \chains^*(F)\langle -p+1,0 \rangle \to 
\cell^R_k \chains^*(F)\langle -p,0 \rangle \Big) \bigg)\]
For $p=0$, we get: $E^1_{0,q}=\pi_q \big(\cell^R_k \chains^*(F)\langle 0,0 \rangle \big)$. 
For $p>0$, note that taking $k$-cellular approximations preserves triangles. Hence:
\begin{align*}
&\cone \big(\cell^R_k \chains^*(F)\langle -p+1,0 \rangle \to \cell^R_k \chains^*(F)\langle -p,0 \rangle \big)
\simeq\\
&\simeq \cell^R_k \cone \big(\chains^*(F)\langle -p+1,0 \rangle \to \chains^*(F)\langle -p,0 \rangle \big)\\
&\simeq \cell^R_k \big( \Sigma^{-p} H^p(F;k) \big)
\end{align*}
The same arguments as in the proof of 
theorem~\ref{the: Convergence of EMSS to fixed points over finite nilpotent space} show that
\[\cell_k^R H^p(F;k) \simeq \cell^N  H^p(F;k)\] for every $p$. This yields the desired $E^1$-page.
\end{proof}
\begin{remark}
The proposition above holds also when $k$ is a commutative, connective,
bounded above $\sphere$-algebra.
\end{remark} 

\subsection{Fibrations over a nilpotent space with cyclic fundamental group}
\label{sub: Fibrations over a nilpotent space with cyclic fundamental group}
The setting we consider is as follows.
Let $F \to E \to B$ be a homotopy fibration sequence. We assume $B$ is a connected, finite, nilpotent space
with a cyclic fundamental group and $E$ is a connected space. Let $k$ be a Noetherian commutative ring.
We use $C$ to denote the fundamental group $\pi_1(B)$.

As mentioned earlier, a result of Dwyer and Greenlees \cite{Dwyer Greenlees} connects $k$-cellular
approximations over $kC$ with $I$-local cohomology of $kC$-modules, where $I$ is the augmentation ideal
of $kC$. This gives the following explicit description for $k$-cellular approximations over $kC$.
Let $x$ be a generator of $C$, so $z=(1-x)$ is a generator of $I$. Then, for every $kC$-chain complex $M$,
there is a triangle:
\[ \cell^C M \to M \to M[1/z]\]
where $M[1/z]$ is the homotopy colimit of the telescope:
$M \xrightarrow{z}M \xrightarrow{z}\cdots$. Note that this implies $(\pi_i M)[1/z]=\pi_i(M[1/z])$.
In particular, we have a long exact sequence:
\[ \cdots \to \pi_i\cell^C M \to \pi_i M \to \pi_i M[1/z] \to \pi_{i-1}\cell^C M \to \cdots \]

We shall use the notation $\Gamma_I M$ for the $I$-power torsion sub-chain complex of $M$, namely:
\[ \Gamma_I M=\{ m \in M| I^n m=0 \text{ for some }n\}\]
So, if $\pi_i(M)=0$ for all $i\neq0$ and $\pi_0(M)=M_0$ we have the following exact sequence:
\[ 0 \to \pi_0\cell^C M \to M_0 \to M_0[1/z] \to \pi_{-1} \cell^C M \to 0\]
This implies 
\[ \pi_0 \cell^C M = \Gamma_I M_0, \quad \pi_{-1} \cell^C M = (M_0[1/z])/M_0 \] 
and $\pi_i \cell^C M=0$ for $i \neq 0,-1$.
We have gathered all the ingredients to prove the following result.
\begin{lemma}
\label{lem: exact sequence for fibration with cyclic pi1}
Let $F \to E \to B$ be a homotopy fibration sequence. Suppose $B$ and $E$ are connected
and $B$ is a finite, nilpotent space with a cyclic fundamental group with generator $x$.
Let $k$ be a Noetherian commutative ring and let $I$ denote the augmentation ideal of the group ring $k[\pi_1(B)]$.
Then:
\[ \pi_0 (\chains^*(E)\otimes_{\chains^*(B)} k) =\Gamma_I H^0(F;k)\]
and for every $n>0$ there is a short exact sequence:
\[ 0 \to H^{n-1}(F;k)[1/z]/H^{n-1}(F;k)
\to \pi_{-n}(\chains^*(E)\otimes_{\chains^*(B)} k) \to \Gamma_I H^{n}(F;k) \to 0\]
where $z=1-x$.
\end{lemma}
\begin{proof}
Since $B$ is a finite space, by lemma~\ref{lem: fibration EMSS gives cellularization} we see that
$\chains^*(E)\otimes_{\chains^*(B)} k \simeq \cell^{k[\kangr B]}_k \chains^*(F)$.
Now we apply the spectral sequence of proposition~\ref{pro: Spectral Sequence}. Note that 
the homotopy groups of $\cell^C H^p(F;k)$ are: $\Gamma_I H^p(F;k)$ in dimension 0, $H^p(F;k)[1/z]/H^p(F;k)$
in dimension -1 and zero elsewhere. Using standard spectral sequence arguments we obtain the desired result.
\end{proof}

Here is an example for the use of this lemma.
\begin{corollary}
\label{cor: when pi1B is cyclic of order 2}
Let $F,E$ and $B$ be as in lemma~\ref{lem: exact sequence for fibration with cyclic pi1}. Suppose in addition
that $\pi_1(B)$ is the cyclic group of order 2 and $F$ is of finite type. If $\chains^*(E;\Int)\otimes_{\chains^*(B;\Int)} \Int \simeq \Int$ then for every $n>0$, $H^n(F;\Int)$ are finite groups of odd order with the action of $\pi_1(B)$ being multiplication by -1, and $H^0(F;\Int)=\Int$ has the trivial action 
of $\pi_1(B)$.
\end{corollary}
\begin{proof}
Denote by $x$ the generator of $\pi_1(B)$ and set $z=1-x \in \Int[\pi_1(B)]$. Since $z^2=2z$, then for every
discrete $\Int[\pi_1(B)]$-module $M$, the discrete module $M[1/z]$ is uniquely divisible by 2.
From lemma~\ref{lem: exact sequence for fibration with cyclic pi1} we see that for $n>0$,
$H^n(F;\Int) \cong H^n(F;\Int)[1/z]$ and in particular $\Gamma_I H^{n}(F;\Int)=0$.
Hence $H^n(F;\Int)$ is a finitely generated abelian group that is 
uniquely divisible by 2, and therefore must have odd order. For any element $a \in H^n(F;\Int)$ the element 
$x\cdot a+a$ is in $\Gamma_I H^{n}(F;\Int)$, therefore $x\cdot a=-a$.

For $H^0(F;\Int)$, note that it is isomorphic to the finitely generated free abelian group 
$\Int\pi_0(F)=\oplus_{a \in \pi_0(F)} \Int$ with the $\pi_1(B)$ action coming from the action of $\pi_1(B)$ on 
$\pi_0(F)$. Suppose the action of $\pi_1(B)$ on $\pi_0 (F)$ had a free orbit. Then
$H^0(F;\Int)$ would have a direct summand isomorphic to $\Int[\pi_1(B)]$.
However, the map $\Int[\pi_1(B)] \to \Int[\pi_1(B)][1/z]$ cannot be a surjection - because $\Int[\pi_1(B)]$ is not a 
2-divisible group. Hence the action of $\pi_1(B)$ on $\pi_0(F)$ is trivial and the result follows.
\end{proof}





\end{document}